\documentclass[11pt, reqno]{amsart}

\usepackage{amsmath,amssymb,amsthm, epsfig}
\usepackage{amsmath}
\usepackage{amssymb}
\usepackage{color}
\usepackage{mathtools}
\usepackage{calrsfs}
\DeclareMathAlphabet{\pazocal}{OMS}{zplm}{m}{n}
\usepackage{epsfig}
\usepackage[mathscr]{eucal}
\usepackage[latin1]{inputenc}
\usepackage{enumerate}

\usepackage{xcolor}
\usepackage{hyperref}
\usepackage[nameinlink]{cleveref}
\hypersetup{
	colorlinks,
	linkcolor={blue!50!black},
	citecolor={blue!50!black},
	urlcolor={blue!50!black}
}

\newtheorem{theorem}{Theorem}
\newtheorem{definition}{Definition}

\newtheorem{lemma}{Lemma}
\newtheorem{proposition}{Proposition}
\newtheorem{corollary}{Corollary}
\newtheorem{remark}{Remark}

\date{}
\numberwithin{equation}{section}
\numberwithin{theorem}{section}
\numberwithin{lemma}{section}
\numberwithin{corollary}{section}
\numberwithin{remark}{section} \numberwithin{proposition}{section}
\numberwithin{definition}{section}

\usepackage[margin=1in]{geometry}

\newcommand{\ssubset}{\subset\joinrel\subset}
\newcommand{\dd}{\mathrm{d}}
\newcommand{\supp}{\operatorname{supp}}
\newcommand{\sign}{\operatorname{sign}}
\newcommand{\eff}{\operatorname{eff}}
\newcommand{\loc}{\operatorname{loc}}
\newcommand{\Id}{\operatorname{Id}}
\newcommand{\AC}{\operatorname{AC}}
\newcommand{\lefthalfcup}{\mathbin{\vrule height 1.6ex depth 0pt width
0.13ex\vrule height 0.13ex depth 0pt width 1.3ex}}
\newcommand{\longrightharpoonup}{\relbar\joinrel\rightharpoonup}
\newcommand{\esssup}{\operatorname{ess}\,\operatorname{sup}}

\begin{document}

\title[Lagrangian structure of relativistic Vlasov systems]{On the Lagrangian structure of transport equations: relativistic Vlasov systems}

\author[H. Borrin]{Henrique Borrin}
\address{Departamento de Matem\'{a}tica Pura e Aplicada, Universidade Federal do Rio Grande do Sul, Porto Alegre - RS, Brazil.}
\email{henrique.borrin@ufrgs.br}

\author[D. Marcon]{Diego Marcon}
\address{Departamento de Matem\'{a}tica Pura e Aplicada, Universidade Federal do Rio Grande do Sul, Porto Alegre - RS, Brazil.}
\email{diego.marcon@ufrgs.br}

\begin{abstract}
We study the Lagrangian structure of relativistic Vlasov systems, such as the relativistic Vlasov-Poisson and the relativistic quasi-eletrostatic limit of Vlasov-Maxwell equations. We show that renormalized solutions of these systems are Lagrangian and that these notions of solution, in fact, coincide. As a consequence, finite-energy solutions are shown to be transported by a global flow. Moreover, we extend the notion of generalized solution for ``effective" densities and we prove its existence. Finally, under a higher integrability assumption of the initial condition, we show that solutions have every energy bounded, even in the gravitational case. These results extend to our setting those obtained by Ambrosio, Colombo, and Figalli \cite{vlasovpoisson} for the Vlasov-Poisson system; here, we analyse relativistic systems and we consider the contribution of the magnetic force into the evolution equation.

\bigskip


\bigskip

\end{abstract}

\maketitle

\section{Introduction}

\subsection{Overview}
In this paper, we are interested in the Lagrangian structure of relativistic Vlasov systems. These systems describe  the evolution of a nonnegative distribution function $f:(0,\infty)\times \mathbb{R}^3\times \mathbb{R}^3\longrightarrow [0,\infty)$ under the action of a self-consistent acceleration:
\begin{equation}\label{prin}
	\begin{cases}
		\partial_t f_t +\hat{v}\cdot \nabla_x f_t+( E_t+\hat{v}\times B_t)\cdot \nabla_v f_t=0 & \text{ in }\quad (0,\infty)\times \mathbb{R}^3\times \mathbb{R}^3;\\[5pt]
		\rho_t(x)=\int_{\mathbb{R}^3} f_t(x,v)\,\dd v,\quad J_t(x)=\int_{\mathbb{R}^3}\hat{v} f_t(x,v)\,\dd v& \text{ in }\quad (0,\infty)\times \mathbb{R}^3;\\[5pt]
		E_t(x)=\sigma_E\int_{\mathbb{R}^3}\rho_t(y)K(x-y)\,\dd y& \text{ in }\quad (0,\infty)\times \mathbb{R}^3;\\[5pt]
		B_t(x)=\sigma_B\int_{\mathbb{R}^3}J_t(y)\times K(x-y)\,\dd y& \text{ in }\quad (0,\infty)\times \mathbb{R}^3.
	\end{cases}
\end{equation} Here, $f_t(x,v)$ denotes the distribution of particles with position $x$ and velocity $v$ at time $t$, $\hat{v}\coloneqq (1+|v|^2)^{-1/2}v$ is the velocity of particles (we assume the speed of light is $c=1$), $\sigma_E \in \{0,\pm 1\}$, $\sigma_B \in \{0,1\}$, and $K: \mathbb{R}^{3} \longrightarrow \mathbb{R}^{3}$ is given by $K(x)=(4\pi)^{-1}x/|x|^3$.

Such systems are very important in mathematical physics and appear in a variety of physical models. Typically, $\rho_{t}$ and $J_{t}$ represent the density of particles and the relativistic particle current density and $E_{t}$ and $B_{t}$ the electric and magnetic fields, respectively. We postpone a derivation of \eqref{prin} and a more complete description of its significance to \Cref{deriv-model}, but summarize what \eqref{prin} models depending on the values of $\sigma_E$ and $\sigma_B$: 

\begin{itemize}
	\item \textbf{Relativistic Vlasov-Poisson equations:} charged particles under a self-consistent electric field or particles under a self-consistent electric and gravitational fields with particle charge $q> q_c$ if $\sigma_E=1$, $\sigma_B=0$; motion of galaxy clusters under a gravitational field or particles under a self-consistent electric and gravitational fields with particle charge $q<q_c$ if $\sigma_E=-1$, $\sigma_B=0$ (see, for instance, \cite[Chapter 5]{smoothed} and references therein);
	\item \textbf{Relativistic Vlasov-Biot-Savart equations\footnote{This terminology, albeit not standard, is in analogy to the Vlasov-Poisson system, since the magnetic field obeys the Biot-Savart law.}:} charged particles under a self-consistent magnetic field; particles under a self-consistent quasi-electrostatic (QES) electromagnetic and gravitational fields with particle charge $q=q_c$ if $\sigma_E= 0$ and $\sigma_B=1$;
	\item \textbf{QES relativistic Vlasov-Maxwell equations:} charged particles under a self-consistent QES electromagnetic field; particles under a self-consistent QES electromagnetic and gravitational fields with particle charge $q>q_c$ if $\sigma_E=\sigma_B=1$;
	\item \textbf{Relativistic gravitational Vlasov-Biot-Savart equations:} charged particles under a self-consistent magnetic and gravitational fields; particles under a self-consistent quasi-magnetostatic (QMS) electromagnetic and gravitational fields with particle charge $q<q_c$ if $\sigma_E=-1$ and $\sigma_B=1$. 
\end{itemize}

Note we allow $\sigma_B=\sigma_E=0$, that is, \eqref{prin} to be the linear transport equation, but its theory is classical and we shall not consider it. Moreover, the fact that the critical charge evolution system coincides with the Vlasov-Biot-Savart system suggests that the displacement current $\partial_t E_t$ behaves like a lower order term; see \eqref{QES} in \Cref{deriv-model}. This is well-known in Electrodynamics \cite{jackson}; Maxwell predicted theoretically as a correction of Amp\`{e}re's law. Nonetheless, we show that it behaves like a lower order term in the magnetic potential energy; see \Cref{magneticpositive} and \Cref{remarkpartialE}.

Concerning the existence of classical solutions of \eqref{prin}, we refer to \cite{15,36,66}, where the existence of local solutions for the relativistic Vlasov-Poisson system is established. As mentioned in \cite[Chapter 5, Section 1.5]{smoothed}, very little is known regarding the existence of global solutions for general initial data. However, existence results can be found, for instance, for spherically and axially symmetric initial data; see \cite{29,34}. In the aforementioned results, it is required higher integrability assumptions and moment conditions  on the initial data. To be more physically relevant, it is desired to avoid such hypotheses even though classical solutions may fail to exist. We thus consider renormalized and generalized solutions, which allow us to establish a Lagrangian structure for the system, global existence results, and (under suitable energy bounds) a global in time maximal regular flow, as we explain in the next section.

\subsection{Main results}

For our purposes, a crucial observation is that \eqref{prin} can be written as
\begin{equation}\label{transport}
	\partial_t f_t + \textbf{b}_t\cdot\nabla_{x,y}f_t=0,
\end{equation} 
where, for each fixed $t>0$, the vector field $\textbf{b}_t: \mathbb{R}^6 \longrightarrow \mathbb{R}^6$ is given by $\textbf{b}_t(x,v)=(\hat{v},E_t+\hat{v}\times B_t)$. Moreover, the vector field is divergence-free, since
\[
\nabla_{x,v}\cdot \textbf{b}_t = \nabla_v\cdot(\hat{v}\times B_t)=(\nabla_v \times\hat{v})\cdot B_t - \hat{v}\cdot (\nabla_v \times B_t)=0.
\]
By the transport nature of \eqref{transport}, it is expected that solutions have a Lagrangian structure, meaning that the initial condition $f_{0}$ is transported to $f_{t}$ by an associated flow. In the weak regularity regime, however, the existence of such flow is not guaranteed by the classical Cauchy-Lipschitz theory. Indeed, since $K$ is locally integrable, we have $E_t,\, B_t \in L^1_{\loc}(\mathbb{R}^3;\mathbb{R}^3)$ whenever $f_t\in L^1(\mathbb{R}^6)$, so that $\textbf{b}_t$ is only in $L^1_{\loc}(\mathbb{R}^6;\mathbb{R}^6)$.

Since $\textbf{b}_t$ is divergence-free, \eqref{transport} can be rewritten as
\[
\partial_t f_t + \nabla_{x,y}\cdot(\textbf{b}_tf_t)=0.
\] The latter can be interpreted in the distributional sense provided $\textbf{b}_tf_t$ is locally integrable which, however, does not follow only from the assumption $f_t\in L^1(\mathbb{R}^6)$. To treat this problem, we introduce a function $\beta \in C^1(\mathbb{R})\cap L^\infty (\mathbb{R})$ such that
\begin{equation}\label{betatransport}
	\partial_t \beta(f_t)+\nabla_{x,v}\cdot(\textbf{b}_t\beta(f_t))=0
\end{equation}
whenever $f_t$ is a smooth solution of \eqref{transport}. Hence, $\textbf{b}_t\beta(f_t)\in L^1_{\loc}$, which leads to the concept of a renormalized solution; see \cite{DiPerna-Lions-1988}.

\begin{definition}[Renormalized solution]\label{renormalized} 
	For a Borel vector field $\textbf{b}\in L_{\loc}^1( [0,T] \times \mathbb{R}^6;\mathbb{R}^6)$, we say a Borel function $f\in L_{\loc}^1([0,T]\times \mathbb{R}^6)$ is a renormalized solution of \eqref{transport} starting from $f_0$ if \eqref{betatransport} holds in the sense of distributions, that is,
	\begin{equation}\label{prinbeta}
		\int_{\mathbb{R}^6}\phi_0(x,v)\beta(f_0(x,v))\,\dd x\,\dd v
		+\int_0^T \!\! \int_{\mathbb{R}^6} \Big[\partial_t \phi_t(x,v)+\nabla_{x,v}\phi_t(x,v)\cdot\textbf{b}_t(x,v)\Big] \beta(f_t(x,v))\,\dd x\,\dd v\,\dd t=0
	\end{equation}
	for all $\phi\in C_c^1([0,T)\times\mathbb{R}^6)$ and $\beta\in C^1(\mathbb{R})\cap L^\infty(\mathbb{R})$.
	
	Moreover, $f\in L^\infty((0,T);L^1(\mathbb{R}^6))$ is called a renormalized solution of \eqref{prin} starting from $f_0$ if, by setting 
	\[
	\rho_t(x) \coloneqq\int_{\mathbb{R}^3} f_t(x,v)\,\dd v,\quad E_t(x)\coloneqq\sigma_E\int_{\mathbb{R}^3}\rho_t(y)\, K(x-y)\,\dd y,
	\]
	\begin{equation}\label{jBbt}
		J_t(x) \coloneqq\int_{\mathbb{R}^3}\hat{v} f_t(x,v)\,\dd v, \quad B_t(x)\coloneqq\sigma_B\int_{\mathbb{R}^3}J_t(y)\times K(x-y)\,\dd y, \quad \text{and}
	\end{equation}
	
	\[
	\textbf{b}_t(x,v) \coloneqq (\hat{v}, E_t(x)+\hat{v}\times B_t(x)),
	\] we have that $f_t$ satisfies \eqref{prinbeta}, for every $\phi\in C_c^1([0,T)\times\mathbb{R}^6)$, with $\textbf{b}_t$ as in \eqref{jBbt}.
\end{definition}
Observe that the integrability assumption $f_t\in L^1(\mathbb{R}^6)$ is used so that $\rho_t,\, J_t$, $E_t$, and $B_t$ are well defined. From now on, we refer to $E_t$ and $B_t$ as the electric and the magnetic fields, respectively, even though $E_t$ may represent a gravitational field as well; see \Cref{deriv-model}.


Our first main result shows that distributional or renormalized solutions of \eqref{prin} are in fact Lagrangian solutions. This gives a characterization of solutions of \eqref{prin}, since Lagrangian solutions are generally stronger than renormalized or distributional solutions.

\begin{theorem}\label{existencesolution} Let $T>0$ and $f$ be a nonnegative function. Assume $f\in L^\infty([0,T);L^1(\mathbb{R}^6))$ is weakly continuous in the sense that
\[
t\longmapsto \int_{\mathbb{R}^6}f_t\,\varphi\,\dd x\,\dd v \ \text{ is continuous for any } \varphi\in C_c(\mathbb{R}^6).
\]
Assume further that:
\begin{enumerate}[$(i)$]
	\item\label{thmitem1} either $f \in L^\infty((0,T);L^\infty(\mathbb{R}^6))$ and $f_t$ is a distributional solution of \eqref{prin} starting from $f_0$; or
	\item\label{thmitem2} $f_t$ is a renormalized solution of \eqref{prin} starting from $f_0$.
\end{enumerate}
Then, $f_t$ is a Lagrangian solution transported by the Maximal Regular Flow $\boldsymbol{X}(t,x)$ associated to $\textbf{b}_t(x,v)=(\hat{v},E_t(x)+\hat{v}\times B_t(x))$ (see \Cref{lagrangian}), starting from $0$. In particular, $f_t$ is renormalized.
\end{theorem}

Next, in \Cref{generalized}, we introduce the concept of generalized solutions, which allows the electromagnetic field to be generated by effective densities $\rho^{\eff}$ and $J^{\eff}$. This may be interpreted as particles vanishing from the phase space but still contributing in the electromagnetic field in the physical space. In fact, generalized solutions are renormalized if the number of particles is conserved in time, as follows from \Cref{jeffequalj}.
This indicates that, should renormalized solutions fail to exist, there must be a loss of mass/charge as $\hat{v}$ approaches the speed of light.

Our second main theorem provides, under minimal assumptions on the initial datum, the global existence of generalized solutions.

\begin{theorem}[Existence of generalized solutions]\label{existencegeneral} 
	Let $f_0\in L^1(\mathbb{R}^6)$ be a nonnegative function. Then there exists a generalized solution $(f_t,\,\rho^{\eff}_t,\,J^{\eff}_t)$ of \eqref{prin} starting from $f_0$ (see \Cref{generalized})). Moreover, the map
\[
t \in [0,\infty) \longmapsto f_t\in L^1_{\loc}(\mathbb{R}^6)
\]
is continuous and the solution $f_t$ is transported by the Maximal Regular Flow associated to field $\boldsymbol{b}^{\eff}_t(x,v)=(\hat{v},E_t^{\eff}+\hat{v}\times B^{\eff}_t)$.
\end{theorem}

In view of \Cref{existencegeneral}, if we assume higher integrability on the initial datum and bounded initial energy, we can prove the existence of a global Lagrangian solution. Moreover, we show strong continuity of densities and fields and that each energy remains bounded in later times. Furthermore, we emphasize that our result holds even in the gravitational case $\sigma_{E} = -1$.

\begin{theorem}[Existence of global Lagrangian solution]\label{finalthm} 
Let $f_0$ be a nonnegative function with every energy bounded (see \Cref{boundedenergy}). Then there exists a global Lagrangian (hence renormalized) solution $f_t\in C([0,\infty);L^1(\mathbb{R}^6))$ of \eqref{prin} with initial datum $f_0$, and the flow is globally defined on $[0,\infty)$ for $f_0$-almost every $(x,v)\in\mathbb{R}^6$, with  $f_t$ being the image of $f_0$ through the incompressible flow.

Moreover, the following properties hold:
\begin{enumerate}[(i)]
	\item the densities $\rho_t,\,J_t$ and the fields $E_t,\, B_t$ are strongly continuous in $L^1_{\loc}(\mathbb{R}^6)$;
	\item for every $t\geq 0$, $f_t$ has every energy bounded independently of time.
\end{enumerate}
\end{theorem}


\subsection{Structure of the paper}

The paper is organized as follows. As previously mentioned, a discussion of the physical interpretation of \eqref{prin} is presented in \Cref{deriv-model}. In \Cref{sec:Lagrangian}, we prove \Cref{existencesolution}. More explicitly, we rely on the machinery for nonsmooth vector fields developed in \cite{existenceflow} to prove the equivalence of renormalized and Lagrangian solutions.
Moreover, in \Cref{maincorollary}, we show that if the electromagnetic and relativistic energies are integrable in $[0,T]$, then its associated flow is globally defined in time.
In \Cref{sec:generalized}, we extend the notion of generalized solutions from \cite[Definition 2.6]{vlasovpoisson} to our setting (see \Cref{generalized}) in order to allow an ``effective" density current of particles (along with the corresponding ``effective" density of particles) and we prove the existence of a Lagrangian solution with the ``effective" acceleration (\Cref{existencegeneral}). Finally, in \Cref{sec:finite-energy}, we prove \Cref{finalthm} under the condition of each bounded energy (see \Cref{boundedenergy}), obtaining a globally defined flow and a solution of \eqref{prin} for all range of $\sigma_E$ and $\sigma_B$.

\subsection*{Aknowledgements} Henrique is partially supported by CAPES through a Master's scholarship. Diego is partially supported by CNPq-Brazil through grant 311354/2019-0.

\section{Lagrangian solution and associated flow}\label{sec:Lagrangian}

In this section, we prove \Cref{existencesolution} which says that Lagrangian and renormalized solutions of \eqref{prin} are equivalent.
For this, we use the machinery developed in \cite[Sections 4 and 5]{vlasovpoisson} combined with a version of \cite[Theorem 4.4]{vlasovpoisson} that we show holds for our vector field $\textbf{b}$ as well.
From now on, we denote by $\mathcal{M}$ the space of measures with finite total mass, by $\mathcal{M}_+$ the space of nonnegative measures with finite total mass, by $\AC(I;\mathbb{R}^6)$ the space of absolutely continuous curves on the interval $I$ with values in $\mathbb{R}^6$, and by $\mathcal{L}^6$ the Lebesgue measure in $\mathbb{R}^6$. We begin with the preliminary definitions of renormalized solutions, and of regular and maximum regular flows:

\begin{definition}[Regular flow]\label{rf} \textnormal{Fix $\tau_1<\tau_2$ and $B \subseteq \mathbb{R}^6$ a Borel set. For a Borel vector field $\textbf{b}:(\tau_1,\tau_2)\times\mathbb{R}^6\longrightarrow \mathbb{R}^6$, we say that $\textbf{X}:[\tau_1,\tau_2]\times B\longrightarrow \mathbb{R}^6$ is a regular flow with vector $\textbf{b}$ when}
\begin{enumerate}[$(i)$]
\item  \textnormal{for a.e. $x\in B$, we have that $\textbf{X}(\cdot,x)\in \AC([\tau_1,\tau_2];\mathbb{R}^6)$ and that it solves the equation $\dot{x}(t)=\textbf{b}_t(x(t))$ a.e. in $(\tau_1,\tau_2)$ with initial condition $\textbf{X}(x,\tau_1)=x$;}
\item \textnormal{there exists $C > 0$ such that $\textbf{X}(t,\cdot)_\#(\mathcal{L}^6\lefthalfcup B)\leq C\mathcal{L}^6$ for all $t\in [\tau_1,\tau_2]$. Note that $C$ can depend on the particular flow $\textbf{X}$.} 
\end{enumerate}
\end{definition}

\begin{definition}[Maximum regular flow]\label{mrf} 
	\textnormal{For every $s\in(0,T)$, a Borel map $\textbf{X}(\cdot,s,\cdot)$ is said to be a maximum regular flow (starting at $s$) if there exist two Borel maps $T^+_{s,\textbf{X}}:\mathbb{R}^6\longrightarrow (s,T]$, $T^-_{s,\textbf{X}}:\mathbb{R}^6\longrightarrow [0,s)$ such that $\textbf{X}(\cdot,s,x)$ is defined in $(T^-_{s,\textbf{X}},T^+_{s,\textbf{X}})$ and}
\begin{enumerate}[$(i)$]
\item \textnormal{for a.e. $x\in\mathbb{R}^6$, we have that $\textbf{X}(\cdot,s,x)\in \AC((T^-_{s,\textbf{X}},T^+_{s,\textbf{X}});\mathbb{R}^6)$ and that it solves the equation $\dot{x}(t)=\textbf{b}_t(x(t))$ a.e. in $(T^-_{s,\textbf{X}},T^+_{s,\textbf{X}})$ with $\textbf{X}(s,s,x)=x$;}
\item \textnormal{there exists a constant $C > 0$ such that $\textbf{X}(t,s,\cdot)_\#(\mathcal{L}^6\lefthalfcup \{T^-_{s,\textbf{X}}<t<T^+_{s,\textbf{X}}\})\leq C\mathcal{L}^6$ for all $t\in [0,T]$. As berfore, this constant $C$ can depend of $\textbf{X}$ and $s$;}
\item \textnormal{for a.e. $x\in\mathbb{R}^6$, either $T^+_{s,\textbf{X}}=T$ and $\textbf{X}(\cdot,s,x)\in C([s,T];\mathbb{R}^6)$, or $\lim_{t\uparrow T^+_{s,\textbf{X}}}|\textbf{X}(t,s,x)|=\infty$. Analogousy, either $T^-_{s,\textbf{X}}=0$ and $\textbf{X}(\cdot,s,x)\in C([0,s];\mathbb{R}^6)$, or $\lim_{t\downarrow T^-_{s,\textbf{X}}}|\textbf{X}(t,s,x)|=\infty$.}
\end{enumerate}
\end{definition}
The following lemma (compare with \cite[Theorem 4.4]{vlasovpoisson}), combined with the facts that $\textbf{b}_t$ is divergence-free in the sense of distribution a.e. in time and that $\boldsymbol{b}\in L^\infty((0,T);L^1_{\loc}(\mathbb{R}^6;\mathbb{R}^6))$, provides a sufficient condition to the existence and the uniqueness of a maximum regular flow for the continuity equation.

\lemma\label{A2} Let 
$\textbf{b}:(0,T)\times \mathbb{R}^6\longrightarrow \mathbb{R}^6
$
 be given by $\textbf{b}_t(x,v)=(\textbf{b}_{1t}(v),\textbf{b}_{2t}(x,v))$, where
\[
\textbf{b}_{1}\in L^\infty((0,T);W^{1,\infty}_{\loc}(\mathbb{R}^3;\mathbb{R}^3)),
\]
\[
\textbf{b}_{2t}(x,v)=K\ast \rho_t(x)+\textbf{b}_{1t}(v)\times\int_{\mathbb{R}^3}K(y-x)\times\dd J_t(y)\eqqcolon K\ast \rho_t(x)+\textbf{b}_{1t}(v)\times \tilde{\textbf{b}}_{2t}(x),
\]
with $\rho\in L^\infty((0,T); \mathcal{M}_+(\mathbb{R}^3))$ and $|J| \in L^\infty((0,T); \mathcal{M}_{+}(\mathbb{R}^3))$.
Then, $\textbf{b}$ satisfies the following: for any nonnegative $\bar{\rho}\in L^\infty(\mathbb{R}^3)$ with compact support and any closed interval $[a,b]\subset [0,T]$, both continuity equations
\[
\frac{\dd}{\dd t}\rho_t\pm\nabla_{x,v}\cdot(\textbf{b}_t\rho_t)=0 \quad \text{in}\, (a,b)\times \mathbb{R}^6
\]
have at most one solution in the class of all weakly* nonnegative continuous functions $[a,b]\ni t\longrightarrow \tilde{\rho}_t$ with $\rho_a=\bar{\rho}$ and $\cup_{t\in[a,b]}\supp \rho_t \Subset \mathbb{R}^6$.
\proof
We proceed as \cite{vlasovpoisson} with the proof for autonomous vector fields (in particular, densities $\rho,\, J$ are independent of time), but the computations generalize for the time dependent case. We denote by $\mathcal{P}(X)$ the set of probability measures on $X$, and by $e_t: C([0,T];\mathbb{R}^6)\longrightarrow \mathbb{R}^6$ the evaluation map at time $t$, which means $e_t(\eta)\coloneqq \eta(t)$. Thanks to \cite[Theorem 5.1]{vlasovpoisson}, any two bounded compactly supported nonnegative distributional solutions with the same initial datum $\bar{\rho}$ can be represented by $\boldsymbol{\eta}_1,\,\boldsymbol{\eta}_2\in \mathcal{P}(C([0,T];B_R\times B_R))$. Setting $\boldsymbol{\eta}=(\boldsymbol{\eta}_1+\boldsymbol{\eta}_2)/2$, if we can prove that the desintegration $\boldsymbol{\eta}_x$ of $\boldsymbol{\eta}$ with respect to $e_0$ is a Dirac delta for $\bar{\rho}$-a.e. $x$ we deduce that $\boldsymbol{\eta}_x=(\boldsymbol{\eta}_1)_x=(\boldsymbol{\eta}_2)_x$ for $\bar{\rho}$-a.e. $x$, which gives $\boldsymbol{\eta}_1=\boldsymbol{\eta}_2$. Hence, it is enough to show that given $\boldsymbol{\eta}\in \mathcal{P}(C([0,T];B_R\times B_R))$ concentrated on integral curves of $\textbf{b}$ such that $(e_t)_\# \boldsymbol{\eta}\leq C_0 \mathcal{L}^{6}$ for all $t\in[0,T]$, $\boldsymbol{\eta}_x$ is  a Dirac delta for $e_{0\#}\boldsymbol{\eta}$-a.e. $x$. For this purpose, consider the function
\[
\Phi_{\delta,\zeta}(t)\coloneqq \iiint \log\left(1+\frac{|\gamma^1(t)-\eta^1(t)|}{\zeta\delta}+\frac{|\gamma^2(t)-\eta^2(t)|}{\delta}\right)\dd\mu(x,\eta,\gamma),
\]
where $\delta,\, \zeta \in(0,1)$ are small constants to be chosen later, $t\in[0,T]$, $\bar{\rho}\coloneqq (e_0)_\#\dd\boldsymbol{\eta}$, $\dd\mu(x,\eta,\gamma)\coloneqq \dd\boldsymbol{\eta}_x(\gamma)\dd\boldsymbol{\eta}_x(\eta)\dd\bar{\rho}(x)$, with notation $\eta(t)=(\eta_1(t),\eta_2(t))\in\mathbb{R}^3\times \mathbb{R}^3$. Note that $\mu\in \mathcal{P}(\mathbb{R}^3\times C([0,T];\mathbb{R}^3)^2)$ and $\Phi_{\delta,\zeta}(0)=0$. We assume by contradition that $\boldsymbol{\eta}_x$ is not a Dirac delta for $\bar{\rho}$-a.e. $x$, which means that there exists a constant $a>0$ such that
\[
\iiint\left(\int_0^T\min\{|\gamma(t)-\eta(t)|,1\}dt\right)\dd\mu(x,\eta,\gamma)\geq a.
\]
Using Fubini's Theorem, the fact that the integrand is bounded by $1$ and $\mu$ has mass $1$, we have that there exist a time $t_0\in(0,T]$ such that
\[
A\coloneqq \left\{(x,\eta,\gamma): \min\{|\gamma(t_0)-\eta(t_0)|,1\}\geq \frac{a}{2T}\right\}
\]
has $\mu$-measure at least $a/(2T)$. Without loss of generality, by assuming $a\leq 2T$, this implies that $|\gamma(t_0)-\eta(t_0)|\geq a/(2T)$ for all $(x,\eta,\gamma)\in A$, hence
\begin{equation}\label{Phibounded}
\Phi_{\delta,\zeta}(t_0)\geq \frac{a}{2T}\log\left(1+\frac{a}{2\delta T}\right).
\end{equation}

Computing the time derivative of $\Phi_{\delta,\zeta}$, we have that
\begin{equation}\label{derivativePhi}
\begin{split}
\frac{\dd\Phi_{\delta,\zeta}}{\dd t}(t)&\leq \iiint\bigg(\frac{|\textbf{b}_1(\gamma^2(t))-\textbf{b}_1(\eta^2(t))|}{\zeta(\delta+|\gamma^2(t)-\eta^2(t)|)}
+\frac{\zeta|\textbf{b}_1(\gamma^2(t))\times(\tilde{\textbf{b}}_2(\gamma^1(t))-\tilde{\textbf{b}}_2(\eta^1(t)))|}{\zeta\delta+|\gamma^1(t)-\eta^1(t)|}
\\
&+\frac{\zeta|(\textbf{b}_1(\gamma^2(t))-\textbf{b}_1(\eta^2(t)))\times \tilde{\textbf{b}}_2(\eta^1(t))|}{\zeta\delta+|\gamma^2(t)-\eta^2(t)|}+\frac{\eta|K\ast\rho(\gamma^1(t))-K\ast\rho(\eta^1(t))|}{\zeta\delta+|\gamma^1(t)-\eta^1(t)|)}\bigg)\dd\mu(x,\eta,\gamma).
\end{split}
\end{equation}
By our assumption on $\textbf{b}_1$, the first summand is easily estimated using the Lipschitz regularity of $\textbf{b}_1$ in $B_R$:
\begin{equation}\label{firstterm}
\iiint\frac{|\textbf{b}_1(\gamma^2(t))-\textbf{b}_1(\eta^2(t))|}{\zeta(\delta+|\gamma^2(t)-\eta^2(t)|)}\dd\mu(x,\eta,\gamma)\leq \frac{\|\nabla\textbf{b}_1\|_{L^\infty(B_R)}}{\zeta}.
\end{equation}
Analogously, the third summand is estimated using the boundedness of $\tilde{\textbf{b}}_2$ and the Lipschitz regularity of $\textbf{b}_1$ in $B_R$:
\begin{equation}\label{thirdterm}
\iiint \frac{\zeta|(\textbf{b}_1(\gamma^2(t))-\textbf{b}_1(\eta^2(t)))\times \tilde{\textbf{b}}_2(\eta^1(t))|}{\zeta\delta+|\gamma^2(t)-\eta^2(t)|} \dd\mu(x,\eta,\gamma)
\leq \zeta\|\nabla\textbf{b}_1\|_{L^\infty(B_R)}\|\tilde{\textbf{b}}_2\|_{L^1(B_R)}.
\end{equation}
For the second term, we have
\[\begin{split}
\iiint \frac{\zeta|\textbf{b}_1(\gamma^2(t))\times(\tilde{\textbf{b}}_2(\gamma^1(t))-\tilde{\textbf{b}}_2(\eta^1(t)))|}{\zeta\delta+|\gamma^1(t)-\eta^1(t)|}&\dd\mu(x,\eta,\gamma)\\
\leq C\|\textbf{b}_1\|_{L^\infty(B_R)}&\iiint \frac{\zeta|K\ast\tilde{\rho}(\gamma^1(t))-K\ast \tilde{\rho}(\eta^1(t))|}{\zeta\delta+|\gamma^1(t)-\eta^1(t)|} \dd\mu(x,\eta,\gamma), 
\end{split}\]
where $\tilde{\rho}(y)\coloneqq \sup_{i}|J_i|(y)$. Since $J_i\in L^\infty((0,\infty); \mathcal{M}(\mathbb{R}^3))$, its total variation is well-defined and has finite measure, thus 
\[
\tilde{\rho}\in L^\infty((0,\infty); \mathcal{M}_+(\mathbb{R}^3)).
\]
By \cite[Theorem 4.4, estimate (4.13)]{vlasovpoisson}, we have that
\begin{equation}\label{secondterm}
\iiint \frac{\zeta|K\ast\bar{\rho}(\gamma^1(t))-K\ast \bar{\rho}(\eta^1(t))|}{\zeta\delta+|\gamma^1(t)-\eta^1(t)|} \dd\mu(x,\eta,\gamma)
\leq C\zeta\left(1+\log\left(\frac{C}{\eta\delta}\right)\right),
\end{equation}
where $\bar{\rho}\in L^\infty((0,\infty); \mathcal{M}_+(\mathbb{R}^3))$ and $C$ depends only on $\bar{\rho}(\mathbb{R}^3)$ and $R$. Hence, the second and fourth terms can be estimated by \eqref{secondterm}.

Then, using \eqref{firstterm}, \eqref{thirdterm}, and \eqref{secondterm}, one can integrate \eqref{derivativePhi} with respect to time in $[0,t_0]$ to obtain
\[
\frac{d\Phi_{\delta,\zeta}}{dt}(t_0)\leq Ct_0\left(\frac{1}{\zeta}+\zeta+\zeta\log\left(\frac{C}{\zeta}\right)+\zeta\log\left(\frac{1}{\delta}\right)\right),
\] 
where $C$ is a constant depending only on $R$, $\rho(\mathbb{R}^3)$, $\tilde{\rho}(\mathbb{R}^3)$, $\|\tilde{\textbf{b}}_2\|_{L^1(B_R)}$, and $\|\textbf{b}_1\|_{W^{1,\infty}(B_R)}$. Choosing first $\zeta>0$ small enough in order to have $Ct_0\zeta< a/(2T)$ and then letting $\delta \longrightarrow 0$, we find a contradiction with \eqref{Phibounded}, concluding the proof.
\endproof
As mentioned before, by \cite[Theorems 5.7, 6.1, 7.1]{existenceflow}, we obtain existence, uniqueness, and a semigroup property for the maximum regular flow (for a concise statement, see \cite[Theorem 4.3]{vlasovpoisson}). We now define generalized flow (analogous to \Cref{rf}) and Lagrangian solutions. For this, we define $\bar{\mathbb{R}}^6=\mathbb{R}^6\cup\{\infty\}$, and given a open set $A\subset [0,\infty)$, $\AC_{\loc}(A;\mathbb{R}^6)$ the set of continuous curves $\gamma:A\longrightarrow \mathbb{R}$ that are absolutely continuous when restricted to any closed interval in $A$.
\begin{definition}[Generalized flow]
For a Borel vector field  $\boldsymbol{b}:(0,T)\times\mathbb{R}^6\longrightarrow \mathbb{R}^6$, the measure $\boldsymbol{\eta}\in\mathcal{M}_+(C[0,T];\bar{\mathbb{R}}^6)$ is said to be a generalized flow of $\boldsymbol{B}$ if $\boldsymbol{\eta}$ is concentrated on the (Borel) set
\[
\Gamma\coloneqq \{\eta\in C([0,T];\bar{\mathbb{R}}^6): \eta\in \AC_{\loc}(\{\eta\neq\infty\};\mathbb{R}^6) \text{ and } \dot{\eta}(t)=\boldsymbol{b}_t(\eta(t)) \text{ for a.e. } t\in \{\eta(t)\neq\infty\}\}.
\]
The generalized flow is regular if there exists $C\geq 0$ such that
\[
(e_t)_\#\boldsymbol{\eta}\lefthalfcup\mathbb{R}^6\leq C \mathcal{L}^6 \quad \forall t\in[0,T].
\]
\end{definition}
\begin{definition}[Transported measures and Lagrangian solutions]\label{lagrangian}
Let $\boldsymbol{b}:(0,T)\times\mathbb{R}^6\longrightarrow \mathbb{R}^6$ be a Borel vector field having a maximal regular flow $\boldsymbol{X}$, and $\boldsymbol{\eta}\in\mathcal{M}_+(C[0,T];\bar{\mathbb{R}}^6)$ with $(e_t)_\#\boldsymbol{\eta}\ll \mathcal{L}^6$ for all $t\in[0,T]$. We say that $\boldsymbol{\eta}$ is transported by $\boldsymbol{X}$ if, for all $s\in[0,T]$, $\boldsymbol{\eta}$ is concentrated on
\[
\{\eta\in C([0,T];\bar{\mathbb{R}}^6):\eta(s)=\infty \text{ or } \eta(\cdot)=\boldsymbol{X}(\cdot,s,\eta(s)) \text{ in } (T^-_{s,\textbf{X}}(\eta(s)),T^+_{s,\textbf{X}}(\eta(s)))\}.
\]
Moreover, let $\rho\in L^\infty((0,T);L^1_{\loc}(\mathbb{R}^6))$ be a nonnegative a distributional solution of the continuity equation, weakly continuous on $[0,T]$ in the duality $C_c(\mathbb{R}^6)$. We say that $\rho_t$ is a Lagrangian solution if there exists $\eta\in\mathcal{M}_+(C([0,T];\bar{\mathbb{R}}^6))$ transported by $\boldsymbol{X}$ with $(e_t)_\#\boldsymbol{\eta}=\rho_t\mathcal{L}^6$ for every $t\in[0,T]$.
\end{definition}
By \cite[Theorem 4.7]{vlasovpoisson}, we have that for $\boldsymbol{b}$ as in \Cref{A2}, regular generalized flows are transported by its maximal regular flow $\boldsymbol{X}$. We are now ready to prove \Cref{existencesolution}.

\proof[Proof of \Cref{existencesolution}]
Notice that the vector field $\boldsymbol{b}$ satisfies $\boldsymbol{b}\in L^\infty((0,T);L^1_{\loc}(\mathbb{R}^6;\mathbb{R}^6))$, is diver-gence-free, and satisfies the uniqueness of bounded compactly supported nonnegative distributional solutions of the continuity equation (see \Cref{A2}). Therefore by \cite[Theorem 5.1]{vlasovpoisson}, we deduce that: if \eqref{thmitem1} holds, then $f_t$ is a Lagrangian solution; if \eqref{thmitem2} holds, then $\beta(f_t)$ is a Lagrangian solution, where $\beta(s)\coloneqq \arctan(s)$. In particular, by \cite[Theorem 4.10]{vlasovpoisson} we have that $f_t$ is a renormalized solution.
\endproof

We have a direct corollary that provides conditions to obtain a globally defined flow, that is, to avoid a finite-time blow up.

\corollary\label{maincorollary} Fix $T>0$ and let $f\in L^\infty((0,T);L^1(\mathbb{R}^6))$ be a nonnegative renormalized solution of \eqref{prin} (as in  \Cref{renormalized}). Assume that
\begin{equation}\label{finiteenergy}
\int_0^T\int_{\mathbb{R}^6}\sqrt{1+|v|^2}f_t(x,v)\,\dd x\,\dd v\,\dd t + \int_0^T  \int_{\mathbb{R}^3} \tfrac{1}{2} |E_t|^2+ \tfrac{1}{2} |B_t|^2\,\dd x\,\dd t<\infty,
\end{equation}
that is, the relativistic energy and the electromagnetic energy \eqref{prin} are integrable in time. 

Then
\begin{enumerate}[(i)]
\item\label{corolitem1} The maximal regular flow $\boldsymbol{X}(t,\cdot)$ associated to $\boldsymbol{b}_t=(\hat{v},E_t+\hat{v}\times B_t)$ and starting from $0$ is globally defined on $[0,T]$ for $f_0$-a.e. $(x,v)$;

\item\label{corolitem2} $f_t$ is the image of $f_0$ through this flow, that is, $f_t=\boldsymbol{X}(t,\cdot)_\#f_0=f_0\circ \boldsymbol{X}^{-1}(t,\cdot)$ for all $t\in[0,T]$:
\[
\int_{\mathbb{R}^6}\phi(x,v)f_t(x,v)\,\dd x\,\dd v=\int_{\mathbb{R}^6}\phi\left(\boldsymbol{X}(t,x,v)\right)f_0(x,v)\,\dd x\,\dd v
\]
for all $ \phi \geq 0, \,  t\in[0,T]$;

\item\label{corolitem3} the map
\[
[0,T]\ni t \longmapsto \int_{\mathbb{R}^6}\psi(f_t(x,v))\,\dd x\,\dd v
\]
is constant in time for all Borel $\psi: [0,\infty)\longrightarrow [0,\infty)$.
\end{enumerate}
\proof
Thanks to \Cref{existencesolution}, the solution is transported by the maximal regular flow associated to $\boldsymbol{b}_t=(v,v\times B_t)$. Moreover, since $f_t$ is renormalized, $g_t\coloneqq \frac{2}{\pi}\arctan f_t : (0,T)\times \mathbb{R}^3\longrightarrow [0,1]$ is a solution of the continuity equation with vector field $\boldsymbol{b}$. Since $g_t^2\leq g_t\leq f_t$ and $|\hat{v}|< 1$, we have
\[\begin{split}
 I&\coloneqq\int_0^T\int_{\mathbb{R}^6}\frac{|\boldsymbol{b}_t(x,v)|g_t(x,v)}{(1+(|x|^2+|v|^2)^{1/2})\log(2+(|x|^2+|v|^2)^{1/2})}\,\dd x\,\dd v\,\dd t\\
&\leq C\int_0^T\int_{\mathbb{R}^6} f_t\, \dd x\,\dd v\,\dd t +\int_0^T\int_{\mathbb{R}^6} \frac{(|E_t|+|B_t|)g_t}{(1+|v|)\log(2+|v|)}\,\dd x\,\dd v\,\dd t\\
&\leq \left(\int_{\mathbb{R}^3}\frac{1}{(1+|v|)^3\log^2(2+|v|)}\dd v\right)\left(\int_0^T\int_{\mathbb{R}^3}|E_t|^2+|B_t|^2\,\dd x\,\dd t\right)+C\int_0^T\int_{\mathbb{R}^6}(1+|v|)f_t\,\dd x\,\dd v\,\dd t.
\end{split}\]
By \eqref{finiteenergy} and $(1+|v|)\leq \sqrt{2(1+|v|^2)}$, we conclude $I$ is bounded.

Now, by the no blow-up criterion in \cite[Proposition 4.11]{vlasovpoisson} we obtain that the maximal regular flow $\boldsymbol{X}$ of $\boldsymbol{b}$ is globally defined on $[0,T]$ (hence, it follows \eqref{corolitem1}). Moreover, the trajectories $\boldsymbol{X}(\cdot,x,v)$ belong to $\AC([0,T];\mathbb{R}^6)$ for $g_0$-a.e. $(x,v)\in \mathbb{R}^6$, and $g_t=\boldsymbol{X}(t,\cdot)_\#g_0=g_0\circ \boldsymbol{X}^{-1}(t,\cdot)$. Since $f_t=\tan \left(\frac{\pi}{2}g_t\right)$ and the map $[0,1) \ni s\longrightarrow \tan\left(\frac{\pi}{2}s\right)\in[0,\infty)$ is a diffeomorphism, we obtain that $f_t=\boldsymbol{X}(t,\cdot)_\#f_0=f_0\circ \boldsymbol{X}^{-1}(t,\cdot)$ (hence, it follows \eqref{corolitem2}). In particular, for all Borel functions $\psi:[0,\infty)\longrightarrow [0,\infty)$ we have
\[
\int_{\mathbb{R}^6}\psi(f_t)\,\dd x\,\dd v=\int_{\mathbb{R}^6}\psi(f_0)\circ \boldsymbol{X}^{-1}(t,\cdot) \,\dd x\,\dd v=\int_{\mathbb{R}^6}\psi(f_0)\,\dd x\,\dd v,
\]
where the second equality follows by the incompressibility of the flow, which gives \eqref{corolitem3}.
\endproof

\begin{remark} \textnormal{As in \cite[Remark 2.4]{vlasovpoisson}, given $0\leq s\leq t \leq T$, it is possible to reconstruct $f_t$ from $f_s$ by using the flow, that is,
$f_t=\boldsymbol{X}(t,s,\cdot)_\#(f_s).$}
\end{remark}
\section{Existence of generalized solution}\label{sec:generalized}
We now introduce the concept of a generalized solution, which allows the electromagnetic field to be generated by effective densities $\rho^{\eff}$ and $J^{\eff}$. We may interpret it as particles vanishing from the phase space but still contributing in the electromagnetic field in the physical space. Thus, it is natural to assume that $\rho^{\eff}_t$ may be larger than $\rho_t$, but it is bounded by the initial particle density $\rho_0$. Moreover, we assume that the particle current density $J^{\eff}_t$ is relativistic and compatible with $\rho^{\eff}_t$, that is, $|J^{\eff}_t|<\rho^{\eff}_t$ and satisfies the continuity equation (see \eqref{conditionj1}, \eqref{conditionj2}, and \eqref{conditionj3} below).
\begin{definition}[Generalized solution]\label{generalized}

Given $\bar{f}\in L^1(\mathbb{R}^6)$, let $f\in L^\infty((0,\infty);L^1(\mathbb{R}^6))$ be a nonnegative function, $\rho^{\eff}_t\in L^\infty((0,\infty); \mathcal{M}_+(\mathbb{R}^3))$, and $(J_t^{\eff})_i\in L^\infty((0,\infty); \mathcal{M}(\mathbb{R}^3))$ for each component $i\in\{1,\,2,\,3\}$.
We say that the triplet $(f_t,\, \rho_t^{\eff},\, J_t^{\eff})$ is a (global in time) generalized solution of \eqref{prin} starting from $\bar{f}$ if, setting
\begin{equation}\label{defeff}
\begin{split}
\rho_t(x)&\coloneqq\int_{\mathbb{R}^3} f_t(x,v)\,\dd v,\quad E^{\eff}_t(x)\coloneqq\sigma_E\int_{\mathbb{R}^3} K(x-y)\,\dd \rho^{\eff}_t(y),\\
J_t(x)&\coloneqq\int_{\mathbb{R}^3}\hat{v} f_t(x,v)\,\dd v, \quad B^{\eff}_t(x)\coloneqq\sigma_B\int_{\mathbb{R}^3} K(y-x)\times\dd J_t^{\eff}(y), \quad \text{and}\\
\textbf{b}^{\eff}_t(x,v)&\coloneqq (\hat{v}, E^{\eff}_t(x)+\hat{v}\times B^{\eff}_t(x)),
\end{split}
\end{equation}
the following hold: $f_t$ is a renormalized solution of the continuity equation with vector field $\boldsymbol{b}_t$ starting from $\bar{f}$,
\begin{subequations}
\begin{equation}\label{conditionj1}
\rho_t\leq \rho^{\eff}_t,\quad |J^{\eff}_t|< \rho^{\eff}_t\quad\text{ as measures for a.e. } t\in(0,\infty),
\end{equation}
\begin{equation}\label{conditionj2}
\rho^{\eff}_t(\mathbb{R}^3)\leq \| f_0\|_{L^1(\mathbb{R}^{6})} \quad \text{ for a.e. } t\in(0,\infty), \text{ and}
\end{equation}
\begin{equation}\label{conditionj3}
\partial_t\rho^{\eff}_t+\nabla\cdot J^{\eff}_t=0\quad \text{with initial condition } \bar{\rho}=\int_{\mathbb{R}^3}\bar{f}\,\dd v\text{, i.e.,}
\end{equation}
\[
\int_{\mathbb{R}^3}\phi_0\,\dd \bar{\rho}+\int_0^\infty\int_{\mathbb{R}^3}(\partial_t\phi_t\,\dd \rho^{\eff}_t+\nabla\phi_t\cdot\dd J^{\eff}_t)\,\dd t=0 \quad \forall \, \phi\in C^1_c([0,\infty)\times\mathbb{R}^3).
\]
\end{subequations}
\end{definition}
Notice that by the Radon-Nikodym's Theorem, combined with \eqref{conditionj1}, there exists a vector field $V^{\eff}\in L^\infty((0,\infty);L^1(\rho^{\eff};\mathbb{R}^3))$ such that $\dd J^{\eff}_t= V^{\eff}_t\,\dd\rho^{\eff}_t$ and $|V^{\eff}_t(x)|< 1$ for a.e. $(t,x)\in(0,\infty)\times\mathbb{R}^3$. This is analogous to the continuity equation associated to \eqref{prin} with initial condition $\rho_0$, which is obtained by integrating \eqref{prin} with respect to $v$ over the whole domain $\mathbb{R}^3$:
\begin{equation}\label{continuityeq}
\int_{\mathbb{R}^3}\phi_0\,\dd \rho_0+\int_0^\infty\int_{\mathbb{R}^3}(\partial_t\phi_t+\nabla\phi_t\cdot V_t)\,\dd \rho_t\,\dd t=0 \quad \forall \phi\in C^1_c([0,\infty)\times\mathbb{R}^3),
\end{equation}
where  $V\coloneqq J/\rho\in L^\infty((0,\infty);L^1(\rho;\mathbb{R}^3))$ satisfies $\dd J_t=V_t\,\dd \rho_t$ and $|V_t(x)|< 1$ for a.e. $(t,x)\in(0,\infty)\times \mathbb{R}^3$. 

To see that \Cref{generalized} is in fact a generalization of \Cref{renormalized}, we remark that $\| \rho_t\|_{L^1(\mathbb{R}^{3})}=\| f_t\|_{L^1(\mathbb{R}^{6})}$, hence it follows by \eqref{conditionj1} and \eqref{conditionj2} that, if the number of particles is conserved a.e. in time, i.e., if $\|f_t\|_{L^1(\mathbb{R}^6)}=\|f_0\|_{L^1(\mathbb{R}^6)}$ for a.e. $t$, then $\rho^{\eff}_t=\rho_t$. Moreover, by \eqref{conditionj3} and \eqref{continuityeq}, we have that $\rho_t$ satisfy the continuity equation with both velocities $V_t$ and $V^{\eff}_t$ with initial condition $\rho_0$. The following lemma gives that $V=V^{\eff}$, whence $J=J^{\eff}$.

\lemma\label{jeffequalj} If $\rho_t$ satisfies the continuity equation with the same initial condition and both vectors $V,\,V^{\eff}$ which
\[
\int_0^T\int_{\mathbb{R}^3}\frac{|\boldsymbol{b}_t(x)|}{1+|x|}\,\dd\mu_t(x)\,\dd t<\infty
\]
holds, then $V=V^{\eff}$.
\proof
Consider a (convex) class $\mathcal{L}_{\boldsymbol{b}}$ of measured-value solutions $\mu_t\in\mathcal{M}_+(\mathbb{R}^3)$ of continuity equation with vector field $\boldsymbol{b}_t$ satisfying
\[
0\leq \partial_t\mu_t\leq \mu_t \quad \Longrightarrow \quad \partial_t\mu_t\in \mathcal{L}_{\boldsymbol{b}}
\]
whenever $\partial_t\mu_t$ still solves the continuity equation with vector field $\boldsymbol{b}_t$, and the integrability condition
\[
\int_0^T\int_{\mathbb{R}^3}\frac{|\boldsymbol{b}_t(x)|}{1+|x|}\,\dd\mu_t(x)\,\dd t<\infty.
\]
Notice that $\rho_t\in \mathcal{L}_{V}\cap \mathcal{L}_{V^{\eff}}$ for all $T>0$, hence by \cite{dipernalions}, we have
\begin{equation}\label{equalityvandveff}
\rho_t=\boldsymbol{X}(t,\cdot)_\#\rho_0=\boldsymbol{X}^{\eff}(t,\cdot)_\#\rho_0 \quad \forall \, t\in[0,T],
\end{equation}
where $\boldsymbol{X}$ and $\boldsymbol{X}^{\eff}$ are $\mathcal{L}_{V}$ and $\mathcal{L}_{V^{\eff}}$ Lagragian flows, respectively, that is, $\boldsymbol{X}(t,\cdot)$ and $\boldsymbol{X}^{\eff}(t,\cdot)$ are (unique) absolutely continuous functions in $[0,T]$ starting from $\rho_0$ (at time $0$) such that 
\[\begin{split}
\dot{\boldsymbol{X}}(t,\cdot)&=V_t(\boldsymbol{X}(t,\cdot)),\quad
\dot{\boldsymbol{X}}^{\eff}(t,\cdot)=V^{\eff}_t(\boldsymbol{X}^{\eff}(t,\cdot)),\\
\boldsymbol{X}(0,\cdot)&=\boldsymbol{X}^{\eff}(0,\cdot)=\Id
\end{split}\]
for $\rho_0$-almost everywhere. By \eqref{equalityvandveff} and the uniqueness of $\boldsymbol{X}$ and $\boldsymbol{X}^{\eff}$, we conclude that $V_t=V^{\eff}_t$.
\endproof

It follows that, if the number of particles is conserved in time, then generalized solutions are renormalized ones. This observation indicates that a generalized solution which is not renormalized must lose mass/charge as the velocity approaches the speed of light.

\medskip

Next, our goal is to prove the global existence of generalized solutions $f_t$ for any nonnegative $f_0\in L^1(\mathbb{R}^6)$ (\Cref{existencegeneral}). In order to do so, we need to establish the existence of a (unique) distributional solution with smooth kernel and initial data. 
More precisely, we show that by smoothing the kernel $K$ and with nonnegative initial condition in $C_c^\infty(\mathbb{R}^6)$, we obtain a classical solution of \eqref{prin}. To avoid any confusion with the notation of \Cref{existencegeneral} and \Cref{finalthm}, we denote by $\pazocal{K}\coloneqq \eta\ast K$ and by $g$ the smoothed kernel and the initial condition, respectively. 

\proposition\label{existencesmooth} Let $g\in C_c^\infty(\mathbb{R}^6)$ be a nonnegative function. Then, there exists a unique nonnegative Lagrangian solution $f\in C^\infty([0,\infty)\times \mathbb{R}^6)$ of the smoothed system \eqref{prin}:
\begin{equation}\label{prinsmooth}
	\begin{cases}
		\partial_t f_t +\hat{v}\cdot \nabla_x f_t+( E_t+\hat{v}\times B_t)\cdot \nabla_v f_t=0 & \text{ in }\quad (0,\infty)\times \mathbb{R}^3\times \mathbb{R}^3;\\
		\rho_t(x)=\int_{\mathbb{R}^3} f_t(x,v)\,\dd v,\quad J_t(x)=\int_{\mathbb{R}^3}\hat{v} f_t(x,v)\,\dd v& \text{ in }\quad (0,\infty)\times \mathbb{R}^3;\\
		E_t(x)=\sigma_E\int_{\mathbb{R}^3}\rho_t(y)\pazocal{K}(x-y)\,\dd y& \text{ in }\quad (0,\infty)\times \mathbb{R}^3;\\
		B_t(x)=\sigma_B\int_{\mathbb{R}^3}J_t(y)\times\pazocal{K}(x-y)\,\dd y& \text{ in }\quad (0,\infty)\times \mathbb{R}^3;\\
		f_0(x,v)=g(x,v)& \text{ in }\quad  \mathbb{R}^3\times \mathbb{R}^3.
	\end{cases}
\end{equation}
\proof 
In this proof, we adapt ideas and techniques from \cite[Chapter 5]{smoothed}.
We construct by induction a sequence of smooth functions $f^n_t$ with initial condition $g$ which converges to a solution of \eqref{prinsmooth}. For $n=1$, let $f^1$ be a solution of the linear transport equation
\[\begin{cases}
	\partial_t f_t^1(x,v)+\nabla_x\cdot(\hat{v} f_t^1)(x,v)=0,\\
	f^1_0(x,v)=g(x,v)
\end{cases}\]
which gives that
\[
f^1_t(x,v)=g(x-t\hat{v},v)\in C_c^\infty([0,\infty)\times\mathbb{R}^6).
\]
Moreover, we have that $f^1$ is a Lagrangian solution, since there exists a unique solution $\boldsymbol{Z}^0(t,\cdot)\coloneqq (\boldsymbol{X}^0,\boldsymbol{V}^0)(t,\cdot)$ of
\[\begin{cases}
	\dot{\boldsymbol{Z}}(t,\cdot)=\boldsymbol{b}^0_t(\boldsymbol{Z}(t,\cdot));\\
	\boldsymbol{Z}(0,\cdot)= \Id,
\end{cases}\]
where $\boldsymbol{b}^0_t(x,v)\coloneqq (\hat{v},0)$. Hence,
\[
f_t^1= g\circ \boldsymbol{Z}^0(t),\quad
\|f_t^1\|_{L^1(\mathbb{R}^6)}=\|g\|_{L^1(\mathbb{R}^6)},\quad\text{and}\quad \|f_t^1\|_{L^\infty(\mathbb{R}^6)}=\|g\|_{L^\infty(\mathbb{R}^6)}.
\]
Now, for $n\geq 2$, assume that there exists a smooth Lagrangian function
\[
f^n\in L^\infty([0,\infty)\times\mathbb{R}^6)\cap\, L^\infty([0,\infty);L^1(\mathbb{R}^6))
\] 
which satisfies
\begin{equation}\label{approx}
	\begin{cases}
		\partial_t f^n_t(x,v)+\nabla_{x,v}\cdot(\boldsymbol{b}^{n-1} f_t^n)(x,v)=0,\\
		f^n_0(x,v)=g(x,v),
	\end{cases}
\end{equation}
where
\[
\boldsymbol{b}^n_t(x,v)=(\hat{v},\,E^n_t+\hat{v}\times B^n_t)(x,v),
\]
and define $f^{n+1}$ as a solution of \eqref{approx} with vector field $\boldsymbol{b}^n_t$. Notice that $\boldsymbol{b}^n_t$ is divergence-free, and since $f^n$ and $\pazocal{K}$ are smooth, we obtain that $\boldsymbol{b}^n_t$ is also smooth. Moreover, we have $\boldsymbol{b}^n\in L^\infty([0,\infty);W^{k,\infty}(\mathbb{R}^6;\mathbb{R}^6))$ for all $k\in \mathbb{N}$, since by Young's inequality (recall that $|J^n|<\rho^n$ a.e.)
\begin{equation}\label{boundbn}
	\begin{split}
		\|D^k_{x,v}\boldsymbol{b}^n_t\|_{L^\infty([0,\infty);L^\infty(\mathbb{R}^6;\mathbb{R}^6))}\leq C\Big(1&+\|K\|_{L^1(B_1;\mathbb{R}^3)}\|D^k \eta\|_{L^\infty(\mathbb{R}^3)}\|\rho^n\|_{L^\infty([0,\infty);L^1(\mathbb{R}^3))}\\
		&+\|K\|_{L^\infty(\mathbb{R}^3\setminus B_1;\mathbb{R}^3)}\|D^k \eta\|_{L^1(\mathbb{R}^3)}\|\rho^n\|_{L^\infty([0,\infty);L^1(\mathbb{R}^3))}\Big).
	\end{split}
\end{equation}
Thus, we have for all $t\geq 0$ a smooth incompressible flow $\boldsymbol{Z}^n(t)=(\boldsymbol{X}^n,\boldsymbol{V}^n)(t)$ which satisfies
\begin{equation}\label{characteristic}
	\begin{cases}
		\dot{\boldsymbol{Z}}(t,\cdot)=\boldsymbol{b}^n_t(\boldsymbol{Z}(t,\cdot));\\
		\boldsymbol{Z}(0,\cdot)= \Id,
	\end{cases}
\end{equation}

and the following properties hold:
\begin{equation}\label{propertiesapprox}
	f_t^{n+1}= g\circ \boldsymbol{Z}^n(t),\quad
	\|f_t^{n+1}\|_{L^1(\mathbb{R}^6)}=\|g\|_{L^1(\mathbb{R}^6)},\quad\text{and}\quad \|f_t^{n+1}\|_{L^\infty(\mathbb{R}^6)}=\|g\|_{L^\infty(\mathbb{R}^6)}.
\end{equation}
Now, we want to exploit the fact that (recall that $g\in C^\infty_c$)
\begin{equation}\label{newidea}
	|f^{n+1}_t-f^n_t|\leq C|\boldsymbol{Z}^n(t)-\boldsymbol{Z}^{n-1}(t)|
\end{equation}
to show that $f^n$ is a Cauchy sequence in $C([0,T]\times\mathbb{R}^6)$. For this purpose, notice that  (we omit the $t$ and $(x,v)$ arguments for a cleaner presentation)
\[\begin{split}
	|\boldsymbol{X}^n(s)-\boldsymbol{X}^{n-1}(s)|&\leq \int_s^t |\boldsymbol{V}^n(\tau)-\boldsymbol{V}^{n-1}(\tau)|+|\boldsymbol{V}^n(\tau)|\left|\frac{1}{\sqrt{1+|\boldsymbol{V}^n(\tau)|^2}}-\frac{1}{\sqrt{1+|\boldsymbol{V}^{n-1}(\tau)|^2}}\right|\,\dd \tau\\
	&\leq \int_s^t |\boldsymbol{V}^n(\tau)-\boldsymbol{V}^{n-1}(\tau)|+\left|\sqrt{1+|\boldsymbol{V}^n(\tau)|^2}-\sqrt{1+|\boldsymbol{V}^{n-1}(\tau)|^2}\right|\,\dd \tau.
\end{split}\]
Thus, by mean value theorem, we conclude that
\[
|\boldsymbol{X}^n(s)-\boldsymbol{X}^{n-1}(s)|\leq 2\int_s^t |\boldsymbol{V}^n(\tau)-\boldsymbol{V}^{n-1}(\tau)|\,\dd \tau.
\]
Moreover, define $E^n$ and $B^n$ as in \eqref{prinsmooth} with densities $\rho^n$ and $J^n$, respectively. Now, by the same procedure as before combined with the uniform boundedness of $B^n$ (by \eqref{boundbn} and \eqref{propertiesapprox}), we have  
\[\begin{split}
	|\boldsymbol{V}^n(s)-\boldsymbol{V}^{n-1}(s)|\leq C\int_s^t |E^n_\tau(\boldsymbol{X}^n(\tau))-E^{n-1}_\tau(\boldsymbol{X}^{n-1}(\tau))|&+|B^n_\tau(\boldsymbol{X}^n(\tau))-B^{n-1}_\tau(\boldsymbol{X}^{n-1}(\tau))|\\
	&+|\boldsymbol{V}^n(\tau)-\boldsymbol{V}^{n-1}(\tau)|\,\dd \tau.
\end{split}\]
By \eqref{boundbn} and \eqref{propertiesapprox}, $E^n$ and $B^n$ are uniformly bounded with respect to $n$ and $t$, thus
\[\begin{split}
	|E^n_\tau(\boldsymbol{X}^n(\tau))-E^{n-1}_\tau(\boldsymbol{X}^{n-1}(\tau))|&\leq |(E^n_\tau-E^{n-1}_\tau)(\boldsymbol{X}^n(\tau))|+|E^{n-1}_\tau(\boldsymbol{X}^n(\tau))-E^{n-1}_\tau(\boldsymbol{X}^{n-1}(\tau))|\\
	&\leq \|E^n_\tau-E^{n-1}_\tau\|_{L^\infty(\mathbb{R}^3)}+C|\boldsymbol{X}^n(\tau)-\boldsymbol{X}^{n-1}(\tau)|,
\end{split}\]
and, analogously,
\[
|B^n_\tau(\boldsymbol{X}^n(\tau))-B^{n-1}_\tau(\boldsymbol{X}^{n-1}(\tau))|\leq \|B^n_\tau-B^{n-1}_\tau\|_{L^\infty(\mathbb{R}^3)}+C|\boldsymbol{X}^n(\tau)-\boldsymbol{X}^{n-1}(\tau)|.
\]
Hence, we obtain that
\[
|\boldsymbol{Z}^n(s)-\boldsymbol{Z}^{n-1}(s)|\leq  C\int_s^t \|E^n_\tau-E^{n-1}_\tau\|_{L^\infty(\mathbb{R}^3)}+ \|B^n_\tau-B^{n-1}_\tau\|_{L^\infty(\mathbb{R}^3)}+|\boldsymbol{Z}^n(\tau)-\boldsymbol{Z}^{n-1}(\tau)|\,\dd \tau.
\]
Thus, by Gronwall's inequality, we conclude that
\[
|\boldsymbol{Z}^n(t)-\boldsymbol{Z}^{n-1}(t)|\leq  C\int_0^t \|E^n_\tau-E^{n-1}_\tau\|_{L^\infty(\mathbb{R}^3)}+ \|B^n_\tau-B^{n-1}_\tau\|_{L^\infty(\mathbb{R}^3)} \,\dd \tau.
\]
Now, by \eqref{propertiesapprox}, we have that $f^n\in C^\infty_c$, which combined with \eqref{newidea} and Young's inequality gives that
\begin{equation}\label{replace}
	\begin{split}
		\|f^{n+1}_t-f^n_t\|_{L^\infty(\mathbb{R}^6)}&\leq C\int_0^t \|\rho^n_\tau-\rho^{n-1}_\tau\|_{L^\infty(\mathbb{R}^3)}+ \|J^n_\tau-J^{n-1}_\tau\|_{L^\infty(\mathbb{R}^3;\mathbb{R}^3)}\,\dd \tau\\
		&\leq C\int_0^t \|f^n_\tau-f^{n-1}_\tau\|_{L^\infty(\mathbb{R}^6)}\,\dd \tau.
	\end{split}
\end{equation}
Therefore, by induction, we have that for all $T>0$,
\[
\|f^{n+1}_t-f^n_t\|_{L^\infty(\mathbb{R}^6)}\leq C\frac{T^n}{n!}, \quad t\in[0,T],
\]
and we conclude that $f^n$ converges uniformly to a function $f\in C([0,\infty)\times\mathbb{R}^6)$. Moreover, by \eqref{propertiesapprox}, we have that $f_t=g\circ \boldsymbol{Z}(t)$, and $f\in L^\infty([0,\infty);L^1(\mathbb{R}^6))\cap L^\infty([0,\infty)\times\mathbb{R}^6)$, where 
\[
\boldsymbol{Z}(t,\cdot)\coloneqq \lim_{n\rightarrow \infty}\boldsymbol{Z}^n(t,\cdot)
\]
Notice that $f_t$ has compact support (since $g\in C^\infty_c$), thus $\rho^n$ and $J^n$ converge to $\rho$ and $J$ in $C([0,\infty)\times\mathbb{R}^6)$, respectively. Therefore, $E^n$ and $B^n$ converge to $E$ and $B$, thus $\boldsymbol{b}^n$ converges to $\boldsymbol{b}$ in $C([0,\infty)\times\mathbb{R}^6)$. By the same computation as \eqref{boundbn}, we have in fact that $\boldsymbol{b}\in C([0,\infty);W^{k,\infty}(\mathbb{R}^6))$ for all $k\in\mathbb{N}$, and we conclude by passing the limit in \eqref{characteristic} that $\boldsymbol{Z}\in C^1([0,\infty); C^\infty(\mathbb{R}^6))$, and we have $f\in C^1([0,\infty); C^\infty(\mathbb{R}^6))$. By iteration, we conclude that $f$ is a smooth nonnegative Lagrangian solution of \eqref{prinsmooth}, where $\boldsymbol{Z}\in C^\infty([0,\infty)\times C^\infty(\mathbb{R}^6))$ solves
\begin{equation}\label{charac}
	\begin{cases}
		\dot{\boldsymbol{Z}}(t,\cdot)=\boldsymbol{b}_t(\boldsymbol{Z}(t,\cdot));\\
		\boldsymbol{Z}(0,\cdot)= \Id.
	\end{cases}
\end{equation}
In particular, we have that $f\in C_c^\infty([0,\infty)\times \mathbb{R}^6)$.

To prove the uniqueness, assume that there exists Lagrangian solutions $f,\,\tilde{f}$ of \eqref{prinsmooth}. Thus, 
\[
f_t\coloneqq g\circ\boldsymbol{Z}(t),\quad \tilde{f}_t\coloneqq g\circ\boldsymbol{\tilde{Z}}(t)
\]
where both $\boldsymbol{Z},\, \boldsymbol{\tilde{Z}}$ solve \eqref{charac}. Thus, we may repeat the proof of \eqref{replace} for $f_t-\tilde{f}_t$ to obtain
\[
\|f_t-\tilde{f}_t\|_{L^\infty(\mathbb{R}^6)}\leq C\int_0^t \|f_\tau-\tilde{f}_\tau\|_{L^\infty(\mathbb{R}^6)}\,\dd \tau,
\]
and we conclude by Gronwall's inequality that $f\equiv\tilde{f}$.
\endproof

We are now able to prove our second main result.

\proof[Proof of \Cref{existencegeneral}]
Our proof follows the same general structure of the proof of \cite[Theorem 2.7]{vlasovpoisson}: we begin by approximating $f$ as a $L^1$ limit of $f^n$ (Steps 1 and 2), which was already shown in \cite{vlasovpoisson}; then, we approximate $(\rho^{\eff}_t,\,J^{\eff}_t)$ and show that the electromagnetic field of the approximation converges to the effective field $(E^{\eff}_t,\,B^{\eff}_t)$ (Steps 3 and 4); finally, in Step 5, we combine stability results for the continuity equation obtained in \cite[Section 5]{vlasovpoisson} to take limits in the approximated system and conclude that the limiting solution is transported by the limit of the incompressible flow.
 
\textbf{Step 1: Approximating solutions.} Consider $K^n\coloneqq K\ast \eta^n$, where $\eta^n(x)\coloneqq n^3\eta(nx)$, and $\eta$ is a standard convolution kernel in $\mathbb{R}^3$. Let $f^n_0\in C^\infty_c(\mathbb{R}^6)$ be a sequence such that 
\begin{equation}\label{convergencel1f0}
f_0^n\longrightarrow f_0 \text{ in } L^1(\mathbb{R}^6).
\end{equation}
Moreover, denote $f^n_t$ the smooth solution of \eqref{prin} with initial condition $f^n_0$ and kernel $K^n$ (see \Cref{existencesmooth}), and its respective charge density, electric field, density current, and magnetic field defined by 
\[\begin{split}
\rho_t^n(x)&\coloneqq \int_{\mathbb{R}^3}f_t^n(x,v)\, \dd v, \quad E^n_t(x)\coloneqq \sigma_E\int_{\mathbb{R}^3} \rho^n_t(y) K^n(x-y)\, \dd y,\\
J^n_t (x)&\coloneqq \int_{\mathbb{R}^3} \hat{v} \,f_t^n(x,v)\, \dd v,\quad \text{and} \quad B^n_t(x)\coloneqq \sigma_B\int_{\mathbb{R}^3} J^n_t(y)\times K^n(x-y)\, \dd y.
\end{split}\]

Since $K^n$ is smooth and vanishes at infinity, we have $E^n,\,B^n\in L^\infty([0,\infty);W^{1,\infty}(\mathbb{R}^3;\mathbb{R}^3))$ (but without a uniform bound with respect to $n$, nonetheless). Hence, $\boldsymbol{b}_t^n\coloneqq (\hat{v},E^n_t+\hat{v}\times B^n_t)$ is a Lipschitz divergence-free vector field, and its flow $\boldsymbol{X}^n(t):\mathbb{R}^6\longrightarrow \mathbb{R}^6$ is well defined and incompressible, hence by theory for the transport equation, for all $t\in [0,\infty)$ and each component $i \in\{1,\, 2,\, 3\}$,
\begin{equation}\label{boundfn}
f_t^n= f_0^n\circ \boldsymbol{X}^n(t)^{-1}\quad\text{ and }\quad
\|J_t^n\|_{L^1(\mathbb{R}^3,\mathbb{R}^3)}\leq\||\hat{v}|f_t^n\|_{L^1(\mathbb{R}^6)}<\|\rho_t^n\|_{L^1(\mathbb{R}^3)}=\|f_t^n\|_{L^1(\mathbb{R}^6)}=\|f_0^n\|_{L^1(\mathbb{R}^6)}.
\end{equation}

Assume without loss of generality that $\mathcal{L}^{6}(\{f_0=k\})=0$ for every $k\in \mathbb{N}$ (otherwise, consider $\mathcal{L}^{6}(\{f_0=k+\tau\})=0$ for $\tau \in (0,1)$), we deduce that for all $k$
\begin{equation}\label{convergencefnk}
f_0^{n,k}\coloneqq \boldsymbol{1}_{\{k\leq f_0^n< k+1\}} f_0^n\longrightarrow f_0^k\coloneqq \boldsymbol{1}_{\{k\leq f_0< k+1\}} f_0 \quad \text{in } L^1(\mathbb{R}^6).
\end{equation}
Thus, by defining $f_t^{n,k}\coloneqq \boldsymbol{1}_{\{k\leq f_t^n< k+1\}} f_t^n$, we have that $f_t^{n,k}$ is a distributional solution of the continuity equation (with vector field $\boldsymbol{b}_t^n$) and $f_0^n$ initial datum. Moreover, we have
\begin{equation}\label{boundjnk}
f_t^{n,k}= \boldsymbol{1}_{\{k\leq f_0^n\circ \boldsymbol{X}^n(t)^{-1}< k+1\}} f_0^n\circ \boldsymbol{X}^n(t)^{-1}, \quad 
\|f_t^{n,k}\|_{L^1(\mathbb{R}^6)}=\|f_0^{n,k}\|_{L^1(\mathbb{R}^6)} 
\quad \forall \, t\in[0,\infty).
\end{equation}

\textbf{Step 2: Limit in phase space.} By construction, $(f^{n,k})_{n\in\mathbb{N}}$ is a nonnegative uniformly bounded sequence. Hence, there exists $f^k\in L^\infty((0,\infty)\times\mathbb{R}^6)$ such that
\begin{equation}\label{weakstarfnk}
f^{n,k}\longrightharpoonup f^k \quad \text{weakly* in } L^\infty((0,\infty)\times\mathbb{R}^6)
\quad\text{as } n\longrightarrow \infty \quad \forall\, k\in\mathbb{N}.
\end{equation}
Moreover, for any $K\ssubset \mathbb{R}^6$, and any bounded function $\phi: (0,\infty)\longrightarrow (0,\infty)$ with compact support, we use test function $\phi(t)\boldsymbol{1}_{K}(x,v)\sign(f_t^k)(x,v)$ for the previous two weak convergence combined with Fatou's lemma, the convergence of $(f^{n,k}_t)_{n\in\mathbb{N}}$, and \eqref{boundjnk} to obtain
\[
\int_0^\infty\phi(t)\|f^{k}_t\|_{L^1(K)}\dd t\leq \left(\int_0^\infty\phi(t)\,\dd t\right)\liminf_{n\rightarrow \infty}\|f^{n,k}_0\|_{L^1(\mathbb{R}^6)}=\left(\int_0^\infty\phi(t)\,\dd t\right)\|f^{k}_0\|_{L^1(\mathbb{R}^6)},
\]
Since $\phi$ was arbitrary the supremum among all compact subset $K\subset \mathbb{R}^6$ we obtain
\begin{equation}\label{boundfk}
\|f_t^k\|_{L^1(\mathbb{R}^6)}\leq\|f_0^k\|_{L^1(\mathbb{R}^6)} \quad \text{for a.e. } t\in(0,\infty), 
\end{equation}
so, in particular, $f^k\in L^\infty((0,\infty);L^1(\mathbb{R}^6))$. Moreover, by defining $f=\sum_{k=0}^\infty f^k$, we have
\begin{equation}\label{regularityf}
\|f_t\|_{L^1(\mathbb{R}^6)}\leq\|f_0\|_{L^1(\mathbb{R}^6)} \quad \text{for a.e. } t\in[0,\infty).
\end{equation}
Noticing that $f^n=\sum_{k=0}^\infty f^{n,k}$, by fixing $\varphi\in L^\infty((0,T)\times\mathbb{R}^6)$, \eqref{boundjnk}, and \eqref{boundfk}, we have for all $k_0\geq 1$,
\[\begin{split}
\left|\int_0^T\int_{\mathbb{R}^6}\varphi(f^n-f)\,\dd x\,\dd v\, \dd t\right|&\leq \left|\sum_{k=0}^{k_0-1}\int_0^T\int_{\mathbb{R}^6}\varphi(f^{n,k}-f^k)\,\dd x\,\dd v\, \dd t\right|\\
&+T\|\varphi\|_{L^\infty((0,T)\times\mathbb{R}^6)}\sum_{k=k_0}^\infty\int_{\mathbb{R}^6}(|f^{n,k}_0|+|f^k_0|)\,\dd x\,\dd v.
\end{split}\]
Now, by the convergence \eqref{weakstarfnk} the first term vanishes as $n\longrightarrow \infty$. Thus, by convergences \eqref{convergencel1f0} and \eqref{convergencefnk}, we have
\[\begin{split}
\limsup_{n\rightarrow \infty}\left|\int_0^T\int_{\mathbb{R}^6}\varphi(f^n-f)\,\dd x\,\dd v\, \dd t\right|\leq 2T\|\varphi\|_{L^\infty((0,T)\times\mathbb{R}^6)}\|f_0\boldsymbol{1}_{\{f_0\geq k_0\}}\|_{L^1(\mathbb{R}^6)}.
\end{split}\]
Letting $k_0\longrightarrow \infty$ and since $\varphi\in L^\infty$ was arbitrary, we conclude
\begin{equation}\label{weakfn}
f^n\longrightharpoonup f \quad \text{ weakly in } L^1((0,T)\times\mathbb{R}^6). 
\end{equation}

\textbf{Step 3: Limit in physical densities.} Since $(\rho^n)_{n\in \mathbb{N}}$ and $(J^n_i)_{n\in \mathbb{N}}$ are bounded sequences in $L^\infty((0,\infty);\mathcal{M}_+(\mathbb{R}^3))$ and $L^\infty((0,\infty);\mathcal{M}(\mathbb{R}^3))$, respectively, for each component $i\in\{1,\,2,\,3\}$ (see \eqref{boundfn}), and $L^\infty((0,\infty);\mathcal{M}(\mathbb{R}^3))= [L^1((0,\infty);C_0(\mathbb{R}^3))]^*$, there exist $\rho^{\eff}\in L^\infty((0,\infty);\mathcal{M}_+(\mathbb{R}^3))$ and $J_i^{\eff}\in L^\infty((0,\infty);\mathcal{M}(\mathbb{R}^3))$ such that
\begin{equation}\label{convergencerhoj}
\begin{split}
\rho^n&\longrightharpoonup \rho^{\eff}\quad \text{weakly* in } L^\infty((0,\infty);\mathcal{M}_+(\mathbb{R}^3));\\
J^n_i&\longrightharpoonup J^{\eff}_i\quad \text{weakly* in } L^\infty((0,\infty);\mathcal{M}(\mathbb{R}^3)).
\end{split}
\end{equation}
for each component $i\in\{1,\,2,\, 3\}$. Hence, by the lower semicontinuity of the norm under weak* convergence, we have
\begin{equation}\label{boundednessrhoeff}
\esssup_{t\in(0,\infty)}|\rho^{\eff}_t|(\mathbb{R}^3)\leq \lim_{n\rightarrow \infty}\left(\sup_{t\in(0,\infty)}\|\rho^n_t\|_{L^1(\mathbb{R}^3)}\right)=\lim_{n\rightarrow \infty}\|\rho^n_0\|_{L^1(\mathbb{R}^3)}=\|f_0\|_{L^1(\mathbb{R}^6)}.
\end{equation}

Now, fixing  a nonnegative function $\varphi\in C_c((0,\infty)\times\mathbb{R}^3)$, by \eqref{weakfn} and \eqref{convergencerhoj}, we obtain that
\[\begin{split}
\int_0^\infty\int_{\mathbb{R}^3}\varphi_t(x)\,\dd \rho_t^{\eff}(x)\,\dd t &\geq\lim_{R\rightarrow\infty}\liminf_{n\rightarrow \infty}\int_0^\infty\int_{\mathbb{R}^3\times B_R} f_t^n(x,v)\varphi_t(x)\,\dd v\,\dd x\,\dd t\\
&=\int_0^\infty\int_{\mathbb{R}^6} f_t(x,v)\varphi_t(x)\,\dd v\,\dd x\,\dd t=\int_0^\infty\int_{\mathbb{R}^3}\varphi_t(x)\,\dd \rho_t(x)\,\dd t.
\end{split}\]
Moreover, by recalling that $|\hat{v}|<1$, we have
\[\begin{split}
\int_0^\infty\int_{\mathbb{R}^3}\varphi_t(x)\,\dd \rho_t^{\eff}(x)\,\dd t &=\lim_{n\rightarrow \infty}\int_0^\infty\int_{\mathbb{R}^6} f_t^n(x,v)\varphi_t(x)\,\dd v\,\dd x\,\dd t\\
&>\lim_{n\rightarrow \infty}\int_0^\infty\int_{\mathbb{R}^6} |\hat{v}|\,f^n_t(x,v)\varphi_t(x)\,\dd v\,\dd x\,\dd t\\
&\geq\int_0^\infty\int_{\mathbb{R}^3}\varphi_t(x)\,\dd |J^{\eff}_t|(x)\,\dd t.
\end{split}\]
Thus,
\begin{equation}\label{boundednessrhoandj}
\rho_t\leq \rho^{\eff}_t,\quad |J^{\eff}_t|< \rho^{\eff}_t\quad\text{ as measures for a.e. } t\in(0,\infty).
\end{equation}

Finally, by the same argument to show \eqref{continuityeq}, we notice that
\[
\int_{\mathbb{R}^3}\phi_0\,\dd \rho^n_0+\int_0^\infty\int_{\mathbb{R}^3}(\partial_t\phi_t\,\dd \rho^n_t+\nabla\phi_t\cdot\,\dd J^n_t)\,\dd t=0 \quad \forall \phi\in C^1_c([0,\infty)\times\mathbb{R}^3).
\]
Hence, by \eqref{convergencel1f0} and \eqref{convergencerhoj}, we conclude by taking the limit $n\longrightarrow \infty$ that
\[
\int_{\mathbb{R}^3}\phi_0\,\dd \rho_0+\int_0^\infty\int_{\mathbb{R}^3}(\partial_t\phi_t\,\dd \rho^{\eff}_t+\nabla\phi_t\cdot\dd J^{\eff}_t)\,\dd t=0 \quad \forall \, \phi\in C^1_c([0,\infty)\times\mathbb{R}^3),
\]
i.e.,
\begin{equation}\label{continuityverified}
\partial_t\rho^{\eff}_t+\nabla\cdot J^{\eff}_t=0 \quad \text{as measures with initial condition } \rho_0.
\end{equation}
\textbf{Step 4: Limit of vector fields.} Using the definition \eqref{defeff}, we claim that
\begin{equation}\label{convergenceb}
\boldsymbol{b}^n\longrightharpoonup \boldsymbol{b}^{\eff} \quad \text{weakly in } L^1_{\loc}((0,\infty)\times\mathbb{R}^6;\mathbb{R}^6)
\end{equation}
and that, for every ball $B_R\subset \mathbb{R}^3$,
\begin{equation}\label{unifomconvergenceb}
[E^n+\hat{v}\times B^n](x+h)\longrightarrow [E^n+\hat{v}\times B^n](x)\text{ as } |h|\rightarrow 0 \text{ in } L^1_{\loc}((0,\infty);L^1(B_R)), \text{ uniformly in } n.
\end{equation}

For this purpose, we first prove that the sequence $(\boldsymbol{b}^n)_{n\in\mathbb{N}}$ is bounded in $L^p_{\loc}((0,\infty)\times\mathbb{R}^6;\mathbb{R}^6)$ for every $p\in[1,3/2)$. Indeed, by using Young's inequality, for every $t\geq 0$, $n\in\mathbb{N}$, and $r>0$,
\[
\|B^n_t\|_{L^p(B_r;\mathbb{R}^3)}+\|E^n_t\|_{L^p(B_r;\mathbb{R}^3)}\leq \|(|J_t^n|\ast \eta^n)\ast |K|\|_{L^p(B_r;\mathbb{R}^3)}+\|(\rho_t^n\ast \eta^n)\ast K\|_{L^p(B_r;\mathbb{R}^3)}
\]
The first term can be bounded by
\[\begin{split}
&\|(|J_t^n|\ast \eta^n)\ast (|K| \boldsymbol{1}_{B_1})\|_{L^p(B_r;\mathbb{R}^3)}+\|(|J_t^n|\ast \eta^n)\ast (|K| \boldsymbol{1}_{\mathbb{R}^3\setminus B_1})\|_{L^p(B_r;\mathbb{R}^3)}\\
&\leq \||J_t^n|\|_{L^1(\mathbb{R}^3)}\|\eta^n\|_{L^1(\mathbb{R}^3)}\|K\|_{L^p(B_1;\mathbb{R}^3)}+\mathcal{L}^3(B_r)^{1/p}\||J_t^n\|_{L^1(\mathbb{R}^3)}\|\eta^n\|_{L^1(\mathbb{R}^3)}\|K\|_{L^\infty(\mathbb{R}^3\setminus B_1;\mathbb{R}^3)}.
\end{split}\]
Likewise, the second term can be bounded by
\[
\|\rho_t^n\|_{L^1(\mathbb{R}^3)}\|\eta^n\|_{L^1(\mathbb{R}^3)}\|K\|_{L^p(B_1;\mathbb{R}^3)}+\mathcal{L}^3(B_r)^{1/p}\|\rho_t^n\|_{L^1(\mathbb{R}^3)}\|\eta^n\|_{L^1(\mathbb{R}^3)}\|K\|_{L^\infty(\mathbb{R}^3\setminus B_1;\mathbb{R}^3)}.
\]
Thus, up to subsequences, the sequence $(\boldsymbol{b}_n)_{n\in\mathbb{N}}$ converges weakly in $L^p_{\loc}$. We now claim that for every $\varphi\in C_c((0,\infty)\times\mathbb{R}^3)$,
\[
\lim_{n\rightarrow \infty}\int_0^\infty\int_{\mathbb{R}^3}(E^n_t+\hat{v}\times B^n_t)\,\varphi_t\,\dd x\,\dd t=\int_0^\infty\int_{\mathbb{R}^3}(E^{\eff}_t+\hat{v}\times B^{\eff}_t)\,\varphi_t\,\dd x\,\dd t.
\]
Indeed, denoting $T_\varphi$ the upper time support of $\varphi$, we have 
\[\begin{split}
&\left|\int_0^\infty\int_{\mathbb{R}^3}(E^n_t+\hat{v}\times B^n_t)\,\varphi_t\,\dd x\,\dd t-\int_0^\infty\int_{\mathbb{R}^3}(E^{\eff}_t+\hat{v}\times B^{\eff}_t)\,\varphi_t\,\dd x\,\dd t\right|\\
&\leq \left|\int_0^\infty\int_{\mathbb{R}^3}(\rho^n_t-\rho^{\eff}_t)\varphi_t\ast K\,\dd x\,\dd t\right|+\left|\int_0^\infty\int_{\mathbb{R}^3}\rho^n_t(\varphi_t\ast K-\varphi_t\ast K\ast \eta^n)\,\dd x\,\dd t\right|\\
&+\left|\int_0^\infty\int_{\mathbb{R}^3}(J^n_t-J^{\eff}_t)\times\varphi_t\ast K\,\dd x\,\dd t\right|+\left|\int_0^\infty\int_{\mathbb{R}^3}J^n_t\times(\varphi_t\ast K-\varphi_t\ast K\ast \eta^n)\,\dd x\,\dd t\right|\\
&\leq \left|\int_0^\infty\int_{\mathbb{R}^3}(\rho^n_t-\rho^{\eff}_t)\varphi_t\ast K\,\dd x\,\dd t\right|+\left|\int_0^\infty\int_{\mathbb{R}^3}(J^n_t-J^{\eff}_t)\times\varphi_t\ast K\,\dd x\,\dd t\right|\\
&+T_\varphi(\|\rho^n\|_{L^\infty((0,\infty);L^1(\mathbb{R}^3))}+\|J^n\|_{L^\infty((0,\infty);L^1(\mathbb{R}^3;\mathbb{R}^3))})\|\varphi\ast K-\varphi\ast K\ast \eta^n\|_{L^\infty((0,\infty)\times\mathbb{R}^3;\mathbb{R}^3)}.
\end{split}\]
By the weak convergence \eqref{convergencerhoj} and the fact that $\varphi\ast K$ is a bounded continuous function, the first and second terms vanish as $n\longrightarrow \infty$. Moreover, the last term also vanishes, since the first factor is bounded by $C\|f_0\|_{L^1(\mathbb{R}^6)}$, where $C>0$ is a universal constant and $\varphi\ast K\ast \eta^n$ convergences uniformly to $\varphi\ast K$ in $(0,\infty)\times \mathbb{R}^3$. Thus, we have proven \eqref{convergenceb}.

We now prove \eqref{unifomconvergenceb}. For this purpose, we combine the fact that $K\in W^{\alpha,p}(\mathbb{R}^3;\mathbb{R}^3)$ for every $\alpha<1$ and $p<3/(2+\alpha)$, and Young's inequality to obtain
\[
\|E_t^n+\hat{v}\times B^n_t\|_{W^{\alpha,p}(B_R;\mathbb{R}^3)}\leq C(R)\|(\rho_t^n+|J^n_t|)\ast\eta^n\|_{L^1(\mathbb{R}^3;\mathbb{R}^3)}.
\]
Combining $\|\eta^n\|_{L^1(\mathbb{R}^3)}=1$ with \eqref{boundfn}, we can bound the right term independently of $n$ and $t$, which combined with the embedding of fractional Sobolev spaces and Nikolsky spaces \cite{nikowsky} gives
\[
\|\boldsymbol{b}^n_t(\cdot+h)-\boldsymbol{b}^n_t(\cdot)\|_{L^p(\mathbb{R}^3;\mathbb{R}^3)}\leq C\left(p,\alpha,R,\|\boldsymbol{b}^n_t\|_{W^{\alpha,p}(B_{2R};\mathbb{R}^3)}\right)|h|^\alpha \quad \forall |h|\leq R,
\]
and \eqref{unifomconvergenceb} follows.

\textbf{Step 5: Conclusion.} By \eqref{convergenceb} and \eqref{unifomconvergenceb}, we can apply the stability result from \cite{dipernalions} to deduce that $f^k$ is a weakly continuous distributional solution of the continuity equation with vector field $\boldsymbol{b}^{\eff}$ and starting from $f_0^k$ for every $k\in\mathbb{N}$. We now exploit the linearity of the continuity equation to show that $F^m\coloneqq \sum_{k=1}^mf^k$ is also a bounded distributional solution for every $m\in \mathbb{N}$. Using the same arguments as in the proof of \Cref{existencesolution}, we obtain that $F^m$ is a renormalized solution for every $m\in\mathbb{N}$. Since $F^m\longrightarrow f$ strongly in $L^1_{\loc}((0,\infty)\times\mathbb{R}^6)$ as $m\longrightarrow\infty$, we obtain that $f$ is a renormalized solution of the continuity equation with vector field $\boldsymbol{b}^{\eff}$ and starting from $f_0$, which combined with \eqref{regularityf} \eqref{boundednessrhoeff}, \eqref{boundednessrhoandj}, and \eqref{continuityverified} proves that the trio $(f_t,\rho^{\eff}_t,J^{\eff}_t)$ is a generalized solution starting from $f_0$ according to \Cref{generalized}.

To show that $f$ is transported by the maximum regular flow associated to $\boldsymbol{b}^{\eff}$, we simply use that each $f^k$ is transported (once again with the same argument as in \Cref{existencesolution}) combined with the definition of $f$ and \eqref{regularityf}. Finally, by \cite[Theorem 4.10]{vlasovpoisson}, we conclude that the map
\[
[0,\infty)\ni t \longmapsto f_t\in L^1_{\loc}(\mathbb{R}^6)\quad \text{is continuous}.\qedhere\]
\endproof

\section{Finite energy solutions}\label{sec:finite-energy}

Up to now, we have established the existence of a generalized solution (see \Cref{existencegeneral}) and that renormalized and generalized solutions coincide in case the mass/charge is conserved in time. In this section, we investigate whether the existence of renormalized solutions can be shown under the more natural condition that the initial total energy is bounded, that is,
\begin{equation}\label{itotalenergy}
\mathcal{E}_0\coloneqq\int_{\mathbb{R}^6}\sqrt{1+|v|^2}f_0(x,v)\,\dd x\,\dd v+\frac{\sigma_E}{2}\int_{\mathbb{R}^3}(H\ast \rho_0) \rho_0 \,\dd x+\frac{\sigma_B}{2}\int_{\mathbb{R}^3}(H\ast J_0) \cdot J_0 \,\dd x<\infty,
\end{equation}
where the first term is the relativistc (initial) total energy and the second and third are the electric and magnetic potential (initial) energies, respectively. For this purpose, we recall that by integrating the first equation of \eqref{prin} with respect to $(x,v)$ on the whole domain $\mathbb{R}^6$ gives that the relativistic energy (formally) satisfies
\[
\frac{\dd}{\dd t}\int_{\mathbb{R}^6}\sqrt{1+|v|^2}f_t(x,v)\,\dd x\,\dd v=\int_{\mathbb{R}^6}\hat{v}\cdot(E_t+\hat{v}\times B_t) f_t(x,v)\,\dd x\,\dd v=\int_{\mathbb{R}^3}E_t\cdot J_t \,\dd x.
\]

Now, Poynting's Theorem gives that the relativistic Vlasov-Maxwell equation has its electromagnetic total energy (formally) conserved, i.e.,
\[
\begin{split}
\int_{\mathbb{R}^6}\sqrt{1+|v|^2}f_t(x,v)\,\dd x\,\dd v+\frac{1}{2}\int_{\mathbb{R}^3}|E_t|^2+|B_t|^2\,\dd x=&\int_{\mathbb{R}^6}\sqrt{1+|v|^2}f_0(x,v)\,\dd x\,\dd v\\
&+\frac{1}{2}\int_{\mathbb{R}^3}|E_0|^2+|B_0|^2\,\dd x,
\end{split}\]
while for the system \eqref{prin} we obtain a similar expression (see \eqref{conservation2} below):
\begin{equation}\label{conservation1}
\begin{split}
\int_{\mathbb{R}^6}\sqrt{1+|v|^2}f_t(x,v)\,\dd x\,\dd v+\frac{\sigma_E}{2}\int_{\mathbb{R}^3}(H\ast \rho_t) \rho_t \,\dd x=&\int_{\mathbb{R}^6}\sqrt{1+|v|^2}f_0(x,v)\,\dd x\,\dd v\\
&+\frac{\sigma_E}{2}\int_{\mathbb{R}^3}(H\ast \rho_0) \rho_0 \,\dd x.
\end{split}
\end{equation}
Notice that the magnetic potential energy does not appear in the conservation above. On the other hand, one can (formally) integrate by parts the electric and magnetic energy to obtain the relations
\begin{equation}\label{conservation2}
\begin{split}
\int_{\mathbb{R}^3}|E_t|^2\,\dd x&=\int_{\mathbb{R}^3}(H\ast \rho_t) \rho_t \,\dd x;\\
\int_{\mathbb{R}^3}|B_t|^2\,\dd x&=\int_{\mathbb{R}^3}(H\ast J_t) \cdot J_t \,\dd x-\int_{\mathbb{R}^3} \left(\nabla\cdot (H\ast J_t)\right)^2\,\dd x,
\end{split}
\end{equation}
where $H(x)\coloneqq (4\pi|x|)^{-1}$. We can interpret $H\ast \rho_t$ and $H\ast J_t$ as the electric potential and magnetic vector potential, respectively (see \cite{jackson}). Notice that, on one hand, the electric potential energy is fully converted into the electric energy. On the other hand, the magnetic potential energy is converted into the magnetic energy and the displacement current $\partial_t E_t$, since
\begin{equation}\label{partialE}
-\int_{\mathbb{R}^3} \left(\nabla\cdot (H\ast J_t)\right)^2\,\dd x=\int_{\mathbb{R}^3} \nabla\cdot (H\ast J_t)\,\partial_t(H\ast \rho_t)\,\dd x=\int_{\mathbb{R}^3} (H\ast J_t)\cdot\partial_tE_t\,\dd x.
\end{equation}

Moreover, we obtain (formally) that the magnetic potential energy is nonnegative for a.e. $t\in[0,\infty)$. Hence, by \eqref{conservation1} and \eqref{conservation2}, we do not expect the initial energy $\mathcal{E}_0$ to bound the total energy of the system given by
\[
\mathcal{E}_t\coloneqq\int_{\mathbb{R}^6}\sqrt{1+|v|^2}f_t(x,v)\,\dd x\,\dd v+\frac{\sigma_E}{2}\int_{\mathbb{R}^3}(H\ast \rho_t) \rho_t \,\dd x+\frac{\sigma_B}{2}\int_{\mathbb{R}^3}(H\ast J_t) \cdot J_t \,\dd x.
\]
Nonetheless, we shall exploit a semicontinuity argument to show an inequality analogous to \eqref{conservation1} (see the proof of \Cref{finalthm}):
\begin{equation}\label{protobound}
\begin{split}
\int_{\mathbb{R}^6}\sqrt{1+|v|^2}f_t(x,v)\,\dd x\,\dd v+\frac{\sigma_E}{2}\int_{\mathbb{R}^3}(H\ast \rho_t) \rho_t \,\dd x \leq & \int_{\mathbb{R}^6}\sqrt{1+|v|^2}f_0(x,v)\,\dd x\,\dd v\\
&+\frac{\sigma_E}{2}\int_{\mathbb{R}^3}(H\ast \rho_0) \rho_0 \,\dd x.
\end{split}
\end{equation}
\begin{remark}\textnormal{Although the formal argument that leads to \eqref{conservation2} suggests the magnetic potential energy is nonnegative, we rigorously justify it in the proof of \Cref{magneticpositive}. Hence, \eqref{itotalenergy} implies that the right-hand side of \eqref{protobound} is bounded.}
\end{remark}
\begin{remark}\label{remarkpartialE}
\textnormal{By \eqref{conservation1} and \eqref{partialE}, we (formally) have}
\begin{equation}\label{lowerorder}
\int_{\mathbb{R}^3}|B_t|^2\,\dd x=\int_{\mathbb{R}^3} A_t \cdot (J_t+\partial_t E_t) \,\dd x,
\end{equation}
\textnormal{where $A_t\coloneqq H\ast J_t$ is the magnetic vector potential. Since we can interpret $\partial_t E_t$ as a density current, one might define the magnetic vector potential as $H\ast (J_t +\partial_t E)$, and therefore \eqref{partialE} does not provide a relation between magnetic energy and magnetic potential energy. We claim that \eqref{lowerorder} still holds if $A_t=H\ast (J_t +\partial_t E)$; thus, we may interpret $\partial_t E_t$ as a lower order term.
Indeed, define a magnetic field with density current $\tilde{J}_t\coloneqq J_t+\partial_t E_t$, that is, $\tilde{B}=\nabla\times(H\ast \tilde{J}_t)$, and a calculation analogous to \eqref{conservation1} gives that
\begin{equation}\label{lowerorder1}
\int_{\mathbb{R}^3}|\tilde{B}_t|^2\,\dd x=\int_{\mathbb{R}^3}(H\ast \tilde{J}_t) \cdot \tilde{J}_t \,\dd x-\int_{\mathbb{R}^3} \left(\nabla\cdot (H\ast \tilde{J})\right)^2\,\dd x.
\end{equation}
Notice that $\nabla\cdot (H\ast \tilde{J})=H\ast(\nabla\cdot J+\partial_t \rho_t)=0$, hence the last term vanishes. Moreover, since $E_t$ is irrotational, $B_t=\tilde{B}_t$; thus, combining \eqref{lowerorder} and \eqref{lowerorder1}, we conclude that
\[
\int_{\mathbb{R}^3}(H\ast \partial_tE_t) \cdot (J_t+\partial_t E_t) \,\dd x=0.
\]
Therefore, had we defined the magnetic vector potential as $H\ast \tilde{J}$, \eqref{lowerorder} would be unaltered.}
\end{remark}

Notice that if $\sigma_E=1$, a bound as \eqref{protobound} gives that each energy term of $\mathcal{E}_t$ is bounded, since $|J|<\rho$ a.e. in space-time. However, it does not provide, in general, control of relativistic energy, electric and magnetic potential energies if $\sigma_E=-1$ or $\sigma_E=0$. If we also assume a higher integrability of $f_0$ and a suitable smallness condition on its norm, the next lemma can be used to bound each energy.

\lemma\label{interpolation} Let $f\in L^1(\mathbb{R}^6)\cap L^q(\mathbb{R}^6)$ be a nonnegative function for some $q\geq 1$ and $\sqrt{1+|v|^2}f\in L^1(\mathbb{R}^6)$. Set $p\coloneqq \frac{4q-3}{3q-2}$. Then $\rho= \int_{\mathbb{R}^3}f(\cdot,v)\,\dd v\in L^p(\mathbb{R}^3)$ and there exists a constant $C>0$, depending only on $q$ such that 
\[
\|\rho\|_{L^p(\mathbb{R}^3)}\leq C\|\sqrt{1+|v|^2}f\|^{\theta}_{L^1(\mathbb{R}^6)}\|f\|^{1-\theta}_{L^q(\mathbb{R}^6)},
\]
where $\theta\coloneqq \frac{3(q-1)}{4q-3}$.
\proof We begin choosing $R>0$ splitting the integral of $\rho$ on the sets $\{|v|<R\}$ and $\{|v|\geq R\}$. Hence, for each $x\in\mathbb{R}^3$,
\[
\rho(x)\leq R^{3(q-1)/q}\|f(x,\cdot)\|_{L^q(\mathbb{R}^3)}+R^{-1}\|\sqrt{1+|v|^2}f(x,\cdot)\|_{L^1(\mathbb{R}^3)}.
\]
By minimizing the right-hand side  with respect to $R$, we have
\[
\rho(x)\leq C\|\sqrt{1+|v|^2}f(x,\cdot)\|^{3(q-1)/(4q-3)}_{L^1(\mathbb{R}^3)}\|f(x,\cdot)\|^{q/(4q-3)}_{L^q(\mathbb{R}^3)}.
\]
Taking the $L^p$-norm on $\rho$ and using H\"{o}lder's inequality, the result follows.
\endproof
As anticipated, if $f_0$ satisfies
\begin{equation}\label{f0bound}
f_0 \in \begin{cases}
L^{1}(\mathbb{R}^6) & \text{if } \sigma_E=1;\\
L^{1}(\mathbb{R}^6)\cap L^{3/2}(\mathbb{R}^6) & \text{if } \sigma_E=0;\\
L^{1}(\mathbb{R}^6)\cap L^{3/2}(\mathbb{R}^6) \text{ and } \|f_0\|_{L^{3/2}(\mathbb{R}^6)}\leq \epsilon & \text{if } \sigma_E=-1
\end{cases}
\end{equation}
for some suitable $\epsilon>0$, the previous lemma allows us to bound  each relativistic energy, electric and magnetic potential energies. Indeed, by Calder\'{o}n-Zygmund estimates and the Sobolev embedding, we have that

\begin{equation}\label{calderonsobolev}
\|H\ast \rho_t\|_{L^6(\mathbb{R}^3)}\leq C\|D^2(H\ast \rho_t)\|_{L^{6/5}(\mathbb{R}^3)}\leq C\|\rho_t\|_{L^{6/5}(\mathbb{R}^3)}
\end{equation}
for some universal constant $C>0$. Combining \eqref{calderonsobolev} with  H\"{o}lder's inequality and \Cref{interpolation} with  $p=6/5$ and $q=3/2$ gives
\begin{equation}\label{potentialenergy}
\begin{split}
\int_{\mathbb{R}^3}(H\ast \rho_t) \rho_t \,\dd x\leq \|H\ast \rho_t\|_{L^6(\mathbb{R}^3)}\|\rho_t\|_{L^{6/5}(\mathbb{R}^3)}&\leq C\|\rho_t\|^2_{L^{6/5}(\mathbb{R}^3)}\\
&\leq C\|\sqrt{1+|v|^2}f_t\|_{L^1(\mathbb{R}^6)}\|f_t\|_{L^{3/2}(\mathbb{R}^6)}.
\end{split}
\end{equation}
Notice that $\|f\|_{L^{\infty}([0,\infty);L^{3/2}(\mathbb{R}^6))}\leq \|f_0\|_{L^{3/2}(\mathbb{R}^6)}$ when the solution is built by approximation (see \eqref{regularityf}). Hence, if \eqref{protobound} holds, we already have a bound of the relativistic energy in the pure magnetic case $\sigma_E=0$, and by the previous bound, we obtain the following boundedness of the magnetic and electric potential energies (recall that $|J|<\rho$ a.e. in space-time):
\[
\int_{\mathbb{R}^3}(H\ast J_t)\cdot J_t \,\dd x\leq  \int_{\mathbb{R}^3}(H\ast \rho_t) \rho_t \,\dd x\leq C\|f_0\|_{L^{3/2}(\mathbb{R}^6)}\int_{\mathbb{R}^6}\sqrt{1+|v|^2}f_0(x,v)\,\dd x\,\dd v.
\]
Now, in the repulsive case $\sigma_E=-1$, we obtain by \eqref{protobound} and \eqref{potentialenergy} that
\[\begin{split}
\left(1-C\|f\|_{L^{\infty}([0,\infty);L^{3/2}(\mathbb{R}^6))}\right)\int_{\mathbb{R}^6}\sqrt{1+|v|^2}f_t(x,v)\,\dd x\,\dd v&\leq \int_{\mathbb{R}^6}\sqrt{1+|v|^2}f_0(x,v)\,\dd x\,\dd v\\
&-\int_{\mathbb{R}^3}(H\ast \rho_0) \rho_0 \,\dd x.
\end{split}\]
Assuming that $f$ is built by approximation as before and that $\|f_0\|_{L^{3/2}(\mathbb{R}^6)}< 1/C\eqqcolon\epsilon$, we have a bound of the relativistic energy; therefore, by \eqref{potentialenergy}, the electric and magnetic potential energies are bounded as well. This motivates the following:
\begin{definition}\label{boundedenergy}
We say that $f_0$ has every energy bounded if \eqref{itotalenergy} and \eqref{f0bound} hold. Moreover, if $f_t$ also satisfies \eqref{protobound} for almost every $t\in[0,\infty)$, then we say that $f_t$ has every energy bounded.
\end{definition}
\begin{remark}
\textnormal{Notice that we need stronger assumptions on the initial data compared to the nonrelativistic Vlasov-Poisson case for $\sigma_E=-1$, where it is only needed that $f_0\in L^{9/7}(\mathbb{R}^3)$, with no smallness assumption (see \cite{thesis}). This is due to the fact that classical kinetic energy grows as $|v|^2$, whereas the relativistic energy as $|v|$.}
\end{remark}
We now prove that if $f_0$ has every energy bounded, then we have a smooth sequence $(f_0^n)_{n\in\mathbb{N}}$ and a mollified sequence of kernels $(H\ast \eta^{k_n})_{n\in\mathbb{N}}$ with uniform bounded energy. We denote by $L_c^\infty$ the space of bounded measurable functions with compact support. 
\lemma\label{approxinitial} Let $\eta^k(x)\coloneqq k^3\eta(kx)$, where $\eta$ is a standard convolution kernel in $\mathbb{R}^3$. Let $f_0$ be a nonnegative function with every energy bounded. Then there exists a sequence $(f_0^n)_{n\in\mathbb{N}}\subset C_c^\infty(\mathbb{R}^6)$ and a sequence $(k_n)_{n\in\mathbb{N}}$ such that $k_n\longrightarrow \infty$ and, by setting $\rho_0^n=\int_{\mathbb{R}^3}f^n_0(\cdot,v)\,\dd v$ and $J_0^n=\int_{\mathbb{R}^3}\hat{v} f^n_0(\cdot,v)\,\dd v$,
\[\begin{split}
\lim_{n\rightarrow \infty}&\left(\int_{\mathbb{R}^6}\sqrt{1+|v|^2}f^n_0(x,v)\,\dd x\,\dd v+\frac{\sigma_E}{2}\int_{\mathbb{R}^3}(H\ast\eta^{k_n}\ast \rho^n_0) \rho^n_0 \,\dd x+\frac{\sigma_B}{2}\int_{\mathbb{R}^3}(H\ast\eta^{k_n}\ast J^n_0) \cdot J^n_0 \,\dd x\right)\\
&=\int_{\mathbb{R}^6}\sqrt{1+|v|^2}f_0(x,v)\,\dd x\,\dd v+\frac{\sigma_E}{2}\int_{\mathbb{R}^3}(H\ast \rho_0) \rho_0 \,\dd x+\frac{\sigma_B}{2}\int_{\mathbb{R}^3}(H\ast J_0) \cdot J_0 \,\dd x.
\end{split}\]
\proof We split the proof in three steps: in Step 1, we assume that $f_0\in L_c^\infty(\mathbb{R}^6)$ and approximate it by a sequence of smooth functions with compact support; in Step 2, we obtain the desired limit without the mollification of $H$; in Step 3, we introduce the mollification of the kernel $\eta^k\ast H$, and conclude that the limit holds if we extract a subsequence of $k$ which depends on $n$.

\textbf{Step 1: $\boldsymbol{f_0\in L_c^\infty(\mathbb{R}^6)}$.} Consider smooth functions $f_0^n$ which converge pointwise such that $\|f_0^n\|_{L^\infty(\mathbb{R}^6)}\leq \|f_0\|_{L^\infty(\mathbb{R}^6)}$ and $\supp f_0^n\subset B_R$ for all $n$ for some $R>0$. Thus, $\|J_0^n\|_{L^\infty(\mathbb{R}^3,\mathbb{R}^3)}<\|\rho_0^n\|_{L^\infty(\mathbb{R}^3)}\leq \|\rho_0\|_{L^\infty(\mathbb{R}^3)}$, and $\supp |J_0^n|\subset\supp\rho_0^n \subseteq B_R$. Moreover, $|H\ast J_0^n|< H\ast \rho_0^n<\infty$ and $H\ast \rho_0^n\longrightarrow H\ast \rho_0$ and $H\ast J_0^n\longrightarrow H\ast J_0$ in $L^p_{\loc}$ for every $p$, and we conclude by dominated convergence that
\begin{equation}\label{approxinitialenergy}
\begin{split}
\lim_{n\rightarrow \infty}&\left(\int_{\mathbb{R}^6}\sqrt{1+|v|^2}f^n_0(x,v)\,\dd x\,\dd v+\frac{\sigma_E}{2}\int_{\mathbb{R}^3}(H\ast \rho^n_0) \rho^n_0 \,\dd x+\frac{\sigma_B}{2}\int_{\mathbb{R}^3}(H\ast J^n_0) \cdot J^n_0 \,\dd x\right)\\
&=\int_{\mathbb{R}^6}\sqrt{1+|v|^2}f_0(x,v)\,\dd x\,\dd v+\frac{\sigma_E}{2}\int_{\mathbb{R}^3}(H\ast \rho_0) \rho_0 \,\dd x+\frac{\sigma_B}{2}\int_{\mathbb{R}^3}(H\ast J_0) \cdot J_0 \,\dd x.
\end{split}
\end{equation}

\textbf{Step 2: $\boldsymbol{f_0\in L^1(\mathbb{R}^6)}$ without mollification of $\boldsymbol{H}$.} By Step 1, it is enough to approximate $f_0$ by $(f_0^n)_{n\in\mathbb{N}}\subset L_c^\infty(\mathbb{R}^6)$ with converging energies to obtain \eqref{approxinitialenergy}. For this purpose, define
\[
f_0^n(x,v)\coloneqq\min\{n,\boldsymbol{1}_{B_n}(x,v)f_0(x,v)\}, \quad (x,v)\in\mathbb{R}^6.
\]
Since $H\geq 0$, the first two integrands on the left-hand side of \eqref{approxinitialenergy} converges monotonically, and we conclude by monotone convergence. Since $|(H\ast J_0^n)\cdot J_0^n|< (H\ast\rho_0) \rho_0$ a.e., and $(H\ast\rho_0)\rho_0$ is integrable (since $f_0$ has every energy bounded), we conclude that the last integral on the left-hand side converges by the dominated convergence.

\textbf{Step 3: Approximation of the kernel.} Given $(f_0^n)_{n\in\mathbb{N}}\in C^\infty_c(\mathbb{R}^6)$ provided by the previous two steps, we have
\[\begin{split}
\lim_{k\rightarrow \infty}&\left(\int_{\mathbb{R}^3}(H\ast\eta^k\ast \rho^n_0) \rho^n_0 \,\dd x+\int_{\mathbb{R}^3}(H\ast\eta^k\ast J^n_0) \cdot J^n_0 \,\dd x\right)\\
&=\int_{\mathbb{R}^3}(H\ast \rho^n_0) \rho^n_0 \,\dd x+\int_{\mathbb{R}^3}(H\ast J^n_0) \cdot J^n_0 \,\dd x
\end{split}\]
for every fixed $n$. Hence, there exists $k_n$ sufficiently large such that
\[\begin{split}
\left|\int_{\mathbb{R}^3}(H\ast\eta^{k_n}\ast \rho^n_0) \rho^n_0 \,\dd x+
\int_{\mathbb{R}^3}(H\ast\eta^{k_n}\ast J^n_0) \cdot J^n_0 \,\dd x-\int_{\mathbb{R}^3}(H\ast \rho^n_0) \rho^n_0 \,\dd x-\int_{\mathbb{R}^3}(H\ast J^n_0) \cdot J^n_0 \,\dd x\right|\\
\leq \frac{1}{n},
\end{split}\]
and the lemma follows.
\endproof

In what follows, we need the following result from \cite[Lemma 3.3]{vlasovpoisson} that we state for convenience of the reader.

\begin{lemma}\label{notproven} Let $T>0$ and $\phi\in C_c((0,T))$ be a nonnegative function. Then, for every sequence $(\rho^n)_{n\in\mathbb{N}}\subset C([0,T];\mathcal{M}_+(\mathbb{R}^3))$ such that 
\[
\sup_{n\in\mathbb{N}}\sup_{t\in[0,T]}\rho^n_t(\mathbb{R}^3)<\infty
\]
and
\begin{equation}\label{rhonconvergence}
\lim_{n\rightarrow \infty}\sup_{t\in[0,T]}\left|\int_{\mathbb{R}^3}\varphi\,\dd (\rho^n_t-\rho_t)\right|=0 \quad \text{for every } \varphi\in C^\infty_c(\mathbb{R}^3).
\end{equation}
we have 
\begin{equation}\label{liminf}
\int_0^T\phi(t)\int_{\mathbb{R}^3}H\ast \rho_t(x)\,\dd\rho_t(x)\,\dd t\leq \liminf_{n\rightarrow \infty}\int_0^T\phi(t)\int_{\mathbb{R}^3}H\ast \eta^{n}\ast \rho^n_t(x)\,\dd\rho^n_t(x)\,\dd t,
\end{equation}
\end{lemma}
Although the previous lemma is enough for $\sigma_E\in\{0,1\}$, we need a slight higher integrability assumption in the gravitational case $\sigma_E=-1$. This is due to the fact that we obtain \eqref{protobound} by a lower semicontinuity argument, and \eqref{liminf} is not sufficient if the electric potential energy is nonpositive. Nonetheless, if $\rho \in L^{6/5}$, we obtain \eqref{liminf} with a limit and an equality, and we prove it in the next lemma.
\lemma\label{proven} Let $\rho^n,\,\rho\in L^\infty([0,T];L^1(\mathbb{R}^3)\cap L^{6/5}(\mathbb{R}^3))$ in the same setting as \Cref{notproven}. Moreover, assume that
\begin{equation}\label{vanishing}
\sup_{n\in\mathbb{N}}\sup_{t\in[0,T]}\|\rho^n_t\|_{L^{6/5}(\mathbb{R}^3)}<\infty.
\end{equation}
Then
\begin{equation}\label{eqproven}
\lim_{n\rightarrow \infty}\int_0^T\phi(t)\int_{\mathbb{R}^3}H\ast \eta^{n}\ast \rho^n_t(x)\,\dd\rho^n_t(x)\,\dd t=\int_0^T\phi(t)\int_{\mathbb{R}^3}H\ast \rho_t(x)\,\dd\rho_t(x)\,\dd t.
\end{equation}
\proof
Notice that
\[\begin{split}
\int_{\mathbb{R}^3}H\ast \eta^n \ast\rho^n_t(x)\rho^n_t(x)-H\ast \rho_t(x)\rho_t(x)\,\dd x&=\int_{\mathbb{R}^3}H\ast \eta^n \ast(\rho^n_t(x)-\rho_t(x))\rho^n_t(x)\,\dd x\\
&+\int_{\mathbb{R}^3}H\ast(\eta^n\ast\rho_t(x)-\rho_t(x))\rho^n_t(x)\,\dd x\\
&+\int_{\mathbb{R}^3}H\ast\rho_t(x)(\rho^n_t(x)-\rho_t(x))\,\dd x\eqqcolon I_1+I_2+I_3.
\end{split}\]
Now, by \eqref{calderonsobolev} and H\"{o}lder inequality, we obtain that 
\[
|I_2|\leq C\|\eta^n \ast\rho_t-\rho_t\|_{L^{6/5}(\mathbb{R}^3)}\sup_{n\in\mathbb{N}}\|\rho^n_t\|_{L^{6/5}(\mathbb{R}^3)}.
\]
Letting $n\longrightarrow \infty$, we obtain that $I_2$ vanishes. We now define $\zeta_k\in C_c^\infty(\mathbb{R}^3)$ as a cutoff function in the annular set $B_{k}\setminus B_{1/k}$, namely,
\[
\begin{cases}
\zeta_k=1 & \text{ in }\quad B_{k}\setminus B_{1/k};\\
\zeta_k=0 & \text{ in }\quad B^c_{k+1}\cup B_{1/(k+1)};\\
0\leq\zeta_k\leq 1& \text{ in }\quad \mathbb{R}^3.
\end{cases}
\]
We write $I_3$ as
\[
|I_3|\leq \left|\int_{\mathbb{R}^3}H\ast\rho_t(x)(\rho^n_t(x)-\rho_t(x))\zeta_k(x)\,\dd x\right|+\left|\int_{\mathbb{R}^3}H\ast\rho_t(x)(\rho^n_t(x)-\rho_t(x))(1-\zeta_k(x))\,\dd x\right|
\]
We want to take first the limit $n\longrightarrow \infty$ and after $k\longrightarrow \infty$ to be able to use \eqref{rhonconvergence}. Now, by \eqref{calderonsobolev}, we obtain
\[\begin{split}
\left|\int_{\mathbb{R}^3}H\ast\rho_t(x)(\rho^n_t(x)-\rho_t(x))(1-\zeta_k(x))\,\dd x\right|&\leq \left|\int_{B_{1/k}\cup B_k^c}H\ast\rho_t(x)(\rho^n_t(x)-\rho_t(x))\,\dd x\right|\\
&\leq  C\|\rho_t\|_{L^{6/5}(\mathbb{R}^3)}\sup_{n\in\mathbb{N}}\|\rho^n_t-\rho_t\|_{L^{6/5}(B_{1/k}\cup B_k^c)}.
\end{split}\]
Defining measures $\dd \mu^n_t\coloneqq (\rho^n_t-\rho_t)^{6/5}\dd x$ and $\mu\coloneqq \sup_{n\in\mathbb{N}}\mu^n$, by \eqref{vanishing} and the continuity from below for measures gives that
\[
\lim_{k\rightarrow \infty}\sup_{n\in\mathbb{N}}\int_{B_{1/k}\cup B_k^c}(\rho^n_t-\rho_t)^{6/5}\,\dd x=\lim_{k\rightarrow \infty}\mu_t(B_{1/k}\cup B_k^c)=\mu_t\left(\cap^\infty_{k=1} B_{1/k}\right)+\mu_t\left(\cap^\infty_{k=1} B^c_k\right)=0,
\]
and we conclude that second term vanishes as $k\longrightarrow \infty$. Now, we bound the first term by
\[
\|H\ast \rho_t\|_{L^\infty(B_{k+1}\setminus B_{1/(k+1)})}\left|\int_{\mathbb{R}^3}\zeta_k(\rho^n_t(x)-\rho_t(x))\,\dd x\right|
\]
By Young's inequality, we have 
\[
\|H\ast \rho_t\|_{L^\infty(B_{k+1}\setminus B_{1/(k+1)})}\leq \|H\|_{L^\infty(B_{k+1}\setminus B_{1/(k+1)})}\|\rho_t\|_{L^1(\mathbb{R}^3)}<\infty.
\]
Hence, by \eqref{rhonconvergence}, $I_3$ vanishes as $n\longrightarrow\infty$ and $k\longrightarrow\infty$. Analogously, we have
\[\begin{split}
|I_1|&=\left|\int_{\mathbb{R}^3}H\ast\eta^n \ast\rho^n_t(x)(\rho^n_t(x)-\rho_t(x))\,\dd x\right|\\
&\leq\|H\ast \rho^n_t\|_{L^\infty(B_{k+1}\setminus B_{1/(k+1)})}\left|\int_{\mathbb{R}^3}\zeta_k(\rho^n_t(x)-\rho_t(x))\,\dd x\right|\\
&+ C\sup_{n\in\mathbb{N}}\|\rho^n_t\|_{L^{6/5}(\mathbb{R}^3)}\sup_{n\in\mathbb{N}}\|\rho^n_t-\rho_t\|_{L^{6/5}(B_{1/k}\cup B_k^c)},
\end{split}\]
and by the same argument as before, $I_1$ vanishes as $n\longrightarrow\infty$ and $k\longrightarrow\infty$, and the lemma follows.
\endproof
We now want to rigorously justify \eqref{conservation2} for $|J| < \rho \in L^{1}(\mathbb{R}^3)$. Actually, the same argument yields the result for $|J| < \rho \in \mathcal{M}_{+}(\mathbb{R}^{3})$. The following lemma gives \eqref{conservation2} with an inequality; in particular, the magnetic potential energy is nonnegative.
\lemma\label{magneticpositive} For every $|J|<\rho\in L^1(\mathbb{R}^3)$ nonnegative,
\begin{equation}\label{conservation2proof}
\begin{split}
\int_{\mathbb{R}^3}|\nabla(H\ast \rho)|^2\,\dd x&\leq\int_{\mathbb{R}^3}(H\ast \rho) \rho \,\dd x;\\
\int_{\mathbb{R}^3}|\nabla\times(H\ast J)|^2\,\dd x&\leq\int_{\mathbb{R}^3}(H\ast J) \cdot J \,\dd x-\int_{\mathbb{R}^3} \left(\nabla\cdot (H\ast J)\right)^2\,\dd x.
\end{split}
\end{equation}
In particular, we obtain that the magnetic potential energy is nonnegative. 
\proof
We split the proof similarly to \Cref{approxinitial}:

\textbf{Step 1: $\boldsymbol{J_i,\,\rho \in L^\infty_c(\mathbb{R}^3)}$.} Consider first $\rho,\, J$ smooth compactly supported functions, and perform an integration by parts to obtain
\[\begin{split}
\int_{B_R}|\nabla(H\ast \rho)|^2\,\dd x&=\int_{B_R}(H\ast \rho) \rho \,\dd x+\int_{\partial B_R}H\ast \rho\, \nabla(H\ast \rho)\cdot \nu_{B_R}\, \dd \mathcal{H}^{2};\\
\int_{B_R}|\nabla\times(H\ast J)|^2\,\dd x&=\int_{B_R}(H\ast J) \cdot J \,\dd x-\int_{\mathbb{R}^3} \left(\nabla\cdot (H\ast J)\right)^2\,\dd x\\
&-\int_{\partial B_R}[(H\ast J) \times (\nabla\times(H\ast J))]\cdot \nu_{B_R}\, \dd \mathcal{H}^{2}\\
&+\int_{\partial B_R}\nabla\cdot(H\ast J)H\ast J\cdot \nu_{B_R}\, \dd \mathcal{H}^{2}.
\end{split}\]
The same identity holds for $J_i,\,\rho \in L^\infty_c(\mathbb{R}^3)$ by approximation for each component $i\in\{1,\,2,\, 3\}$. Since $H\ast \mu$ and $\nabla (H\ast \mu)$ decay as $R^{-1}$ and $R^{-2}$ when evaluated at $\partial B_R$ for all $\mu\in L^\infty_c(\mathbb{R}^3)$, the boundary terms vanish as $R\longrightarrow \infty$, and we obtain that \eqref{conservation2proof} holds with an equality.

\textbf{Step 2: $\boldsymbol{J_i,\,\rho \in L^1(\mathbb{R}^3)}$.} We consider the truncations 
\[\rho^n\coloneqq \min\{n,\boldsymbol{1}_{B_n}(x,v)\rho\}, \quad J_i^n\coloneqq \min\{n,\boldsymbol{1}_{B_n}(x,v)J_i\}.\]
Since $H\geq 0$, by monotone convergence and Step 1 we obtain that
\[
\int_{\mathbb{R}^3}(H\ast \rho) \rho \,\dd x=\lim_{n\rightarrow\infty}\int_{\mathbb{R}^3}(H\ast \rho^n) \rho^n \,\dd x= \lim_{n\rightarrow\infty} \int_{\mathbb{R}^3}|\nabla(H\ast \rho^n)|^2\,\dd x.
\]
Moreover, since $|J|<\rho$, by dominated convergence and Step 1 we obtain that
\[\begin{split}
\int_{\mathbb{R}^3}(H\ast J) \cdot J \,\dd x&=\lim_{n\rightarrow\infty}\int_{\mathbb{R}^3}(H\ast J^n) \cdot J^n \,\dd x\\
&=\lim_{n\rightarrow\infty}\int_{\mathbb{R}^3}|\nabla\times(H\ast J^n)|^2\,\dd x+\lim_{n\rightarrow\infty}\int_{\mathbb{R}^3} \left(\nabla\cdot (H\ast J^n)\right)^2\,\dd x.
\end{split}\]
Assuming without loss of generality that $(H\ast \rho) \rho\in L^1(\mathbb{R}^3)$, we get bounded sequences $(\nabla(H\ast \rho^n))_{n\in\mathbb{N}}$, $(\nabla\cdot(H\ast J^n))_{n\in\mathbb{N}}$, and $(\nabla\times(H\ast J^n))_{n\in\mathbb{N}}$ in $L^2$. Since each sequence converges in the sense of distributions to $\nabla(H\ast \rho)$, $\nabla\cdot(H\ast J)$, and $\nabla\times(H\ast J)$, respectively, and the lower semicontinuity of the $L^2$-norm with respect to the weak convergence, we conclude \eqref{conservation2proof}.
\endproof

Finally, we prove our third main result.

\proof[Proof of \Cref{finalthm}]
The proof of existence of renormalized solutions begins similarly to the proof of \Cref{existencegeneral}: let $(f^n_0)_{n\in\mathbb{N}}\subset C_c^\infty(\mathbb{R}^6)$ and $(k_n)_{n\in\mathbb{N}}$ given by \Cref{approxinitial}. By Steps 1-3 in the proof of \Cref{existencegeneral} we get a sequence of smooth functions $f^n$  satisfying \eqref{prin} with initial condition $f^n_0$ and kernel $K^n$ (see \Cref{existencesmooth}) such that
\begin{equation}\label{convergences}
\begin{split}
f^n\longrightharpoonup f \quad& \text{weakly in } L^1([0,T]\times\mathbb{R}^6) \text{ for any } T>0;\\ 
\rho^n\longrightharpoonup \rho^{\eff} \quad& \text{weakly* in } L^\infty((0,\infty);\mathcal{M}_+(\mathbb{R}^3));\\
|J^{\eff}|<\rho^{\eff} \quad& \text{as measures};\\
\partial_t \rho^{\eff}+\nabla\cdot J^{\eff}=0 \quad& \text{ as measures with initial condition } \rho_0,
\end{split}
\end{equation}
where $\rho^n_t(x)\coloneqq \int_{\mathbb{R}^3}f^n_t(x,v)\,\dd v$. Analogously to \eqref{regularityf}, we have that for $\sigma_E\in\{-1,\, 0\}$,
\begin{equation}\label{ftregularity}
\|f^n_t\|_{L^{3/2}(\mathbb{R}^6)}\leq\|f_0\|_{L^{3/2}(\mathbb{R}^6)},\quad \|f_t\|_{L^{3/2}(\mathbb{R}^6)}\leq\|f_0\|_{L^{3/2}(\mathbb{R}^6)} \quad \text{for a.e. } t\in[0,\infty).
\end{equation}
Moreover, since \eqref{conservation1} holds for classical solutions and $f_0$ has every energy bounded, we obtain that
\begin{equation}\label{boundnenergy}
\sup_{n\in\mathbb{N}}\sup_{t\in[0,\infty)}\int_{\mathbb{R}^6}\sqrt{1+|v|^2}f_t^n\,\dd x\,\dd v\leq C,
\end{equation}
and by the lower semicontinuity of the relativistic energy we deduce that, for every $T>0$,
\begin{equation}\label{boundenergyT}
\int_0^T\int_{\mathbb{R}^6}\sqrt{1+|v|^2}f_t\,\dd x\,\dd v\,\dd t\leq  \liminf_{n\rightarrow \infty}\int_0^T\int_{\mathbb{R}^6}\sqrt{1+|v|^2}f_t^n\,\dd x\,\dd v\,\dd t\leq C T.
\end{equation}
We now claim that $\rho^{\eff}=\rho$ and, consequently, $J=J^{\eff}$, where $|J|<\rho \in L^{\infty}((0,T);L^1(\mathbb{R}^6))$ as in \eqref{prin}. For this, consider $\zeta_k:\mathbb{R}^6\longrightarrow [0,1]$ a nonnegative function which equals $1$ inside $B_k$ and $0$ in $B_{k+1}^c$ and compute
\[\begin{split}
\int_0^\infty\int_{\mathbb{R}^3}(\rho_t^n-\rho_t)\varphi_t\,\dd x\,\dd t&=\int_0^\infty\int_{\mathbb{R}^6}(f_t^n(x,v)-f_t(x,v))\varphi_t(x)\zeta_k(v)\,\dd v\,\dd x\,\dd t\\
&+\int_0^\infty\int_{\mathbb{R}^6} f_t^n(x,v)\varphi_t(x)(1-\zeta_k(v))\,\dd v\,\dd x\,\dd t\\
&+\int_0^\infty\int_{\mathbb{R}^6} f_t(x,v)\varphi_t(x)(\zeta_k(v)-1)\,\dd v\,\dd x\,\dd t.
\end{split}\]
By the weak convergence in $L^1$ in \eqref{convergences}, the first term vanishes as $n\longrightarrow \infty$. The second and third terms can be estimated using \eqref{boundnenergy} and \eqref{boundenergyT}:
\[\begin{split}
&\left|\int_0^\infty\int_{\mathbb{R}^6} f_t^n(x,v)\varphi_t(x)(1-\zeta_k(v))\,\dd v\,\dd x\,\dd t
+\int_0^\infty\int_{\mathbb{R}^6} f_t(x,v)\varphi_t(x)(\zeta_k(v)-1)\,\dd v\,\dd x\,\dd t\right|\\
&+\frac{\|\varphi\|_{L^\infty((0,\infty)\times\mathbb{R}^3)}}{k}\int_0^{T_\varphi}\int_{\mathbb{R}^6}\sqrt{1+|v|^2}f^n_t\,\dd x\,\dd v\,\dd t+\frac{\|\varphi\|_{L^\infty((0,\infty)\times\mathbb{R}^3)}}{k}\int_0^{T_\varphi}\int_{\mathbb{R}^6}\sqrt{1+|v|^2}f_t\,\dd x\,\dd v\,\dd t\\
&\leq \frac{C {T_\varphi}\|\varphi\|_{L^\infty((0,\infty)\times\mathbb{R}^3)}}{k},
\end{split}\]
where $T_\varphi$ is the time support of $\varphi$. Letting $k\longrightarrow \infty$, we conclude that $\rho^n$ converges to $\rho$ weakly* in $L^\infty((0,\infty);\mathcal{M}_+(\mathbb{R}^3))$, which combined with \eqref{convergences} gives that $\rho=\rho^{\eff}$. Hence, by \eqref{convergences} and \Cref{jeffequalj}, we conclude that $J=J^{\eff}$, and in Steps 4 and 5 in the proof of \Cref{existencegeneral}, we obtain a global Lagrangian (hence renormalized) solution $f_t\in C([0,\infty);L^1_{\loc}(\mathbb{R}^6))$ of \eqref{prin} with initial datum $f_0$.

We now prove properties by a lower semicontinuous argument on the energy of $f^n$.

\textbf{Step 1: Bound on the total energy for $\boldsymbol{\mathcal{L}^1}$-almost every time.} We use the weak convergence of $f^n$ (see \eqref{convergences}) with test function $\phi(t)\sqrt{1+|v|^2}\chi_r(x,v)$, where $\phi\in C_c^{\infty}((0,\infty))$ and $\chi_r\in C_c^\infty(\mathbb{R}^6)$ are nonnegative functions, with $\chi_r$ being a cutoff between $B_r$ and $B_{r+1}$, we obtain
\[
\int_0^\infty\int_{\mathbb{R}^6} f_t(x,v)\sqrt{1+|v|^2}\phi(t)\chi_r(x,v)\,\dd v\,\dd x\,\dd t\leq\liminf_{n\rightarrow \infty} \int_0^\infty \phi(t)\int_{\mathbb{R}^6} \sqrt{1+|v|^2}f^n_t(x,v)\,\dd v\,\dd x\,\dd t.
\]
Taking the supremum with respect to $r$, we deduce that
\begin{equation}\label{liminfrelativistic}
\int_0^\infty \phi(t)\int_{\mathbb{R}^6}\sqrt{1+|v|^2} f_t(x,v)\,\dd v\,\dd x\,\dd t\leq\liminf_{n\rightarrow \infty} \int_0^\infty \phi(t)\int_{\mathbb{R}^6} \sqrt{1+|v|^2}f^n_t(x,v)\,\dd v\,\dd x\,\dd t.
\end{equation}
Since $\phi$ is arbitrary, we have that $\sqrt{1+|v|^2}f_t\in L^1_{\loc}(\mathbb{R}^6)$ for almost every $t$. Moreover, since we can decompose the density current as $J=V\rho$ (see remark after \Cref{generalized}), where $|V|<1$ a.e. in spacetime, we have that
\[
\sup_{t\in[0,\infty)}\int_{\mathbb{R}^3}|V_t(x)|\,\dd\rho_t(x)<\infty,
\]
hence by \cite[Theorem 8.1.2]{gradientflows}, we have that $\rho_t$ has a weakly* continuous representative. Furthermore, since $\rho^n$ satisfies a similar continuity equation, by the proof of \cite[Theorem 8.1.2]{gradientflows}, we have that
\[
\left|\int_{\mathbb{R}^3}(\rho_t^n-\rho_s^n)\varphi\,\dd x\right|\leq \|\varphi\|_{C^1(\mathbb{R}^3)}\int_s^t\int_{\mathbb{R}^3}|V^n_r|\rho_r^n\,\dd x\,\dd r\leq C|t-s|
\]
for all $\varphi\in C^\infty_c(\mathbb{R}^3)$, which gives that the map $t\longmapsto\int_{\mathbb{R}^3}\varphi\,\dd \rho^n_t$ is equicontinuous. By the weak* convergence of $\rho^n$ to $\rho$ in $L^\infty((0,\infty);\mathcal{M}_+(\mathbb{R}^3))$, we have a uniform boundedness, thus Arzel\`{a}-Ascoli theorem implies that
\begin{equation}\label{rhonweakconvergence}
\lim_{n\rightarrow \infty}\sup_{t\in[0,T]}\left|\int_{\mathbb{R}^3}\varphi\,\dd (\rho^n_t-\rho_t)\right|=0 \quad \text{for every } \varphi\in C^\infty_c(\mathbb{R}^3).
\end{equation}
Combining the above with the fact that $\rho^n_t$ is uniformly bounded with respect to $n$ and $t$, by \Cref{notproven} we obtain
\begin{equation}\label{liminfelectric}
\int_0^\infty\phi(t)\int_{\mathbb{R}^3}H\ast \rho_t(x)\,\dd\rho_t(x)\,\dd t\leq \liminf_{n\rightarrow \infty}\int_0^\infty\phi(t)\int_{\mathbb{R}^3}H\ast \eta^{k_n}\ast \rho^n_t(x)\,\dd\rho^n_t(x)\,\dd t.
\end{equation}
Combining \eqref{liminfrelativistic}, \eqref{liminfelectric}, and \eqref{conservation1}, we conclude that for $\sigma_E\in\{0,\,1\}$
\[\begin{split}
&\int_0^\infty\phi(t)\left(\int_{\mathbb{R}^6}\sqrt{1+|v|^2} f_t(x,v)\,\dd v\,\dd x+\frac{\sigma_E}{2}\int_{\mathbb{R}^3}H\ast \rho_t(x)\rho_t(x)\,\dd x\right)\,\dd t\\
&\leq\liminf_{n\rightarrow \infty}\int_0^\infty\phi(t)\left(\int_{\mathbb{R}^6}\sqrt{1+|v|^2} f^n_0(x,v)\,\dd v\,\dd x+\frac{\sigma_E}{2}\int_{\mathbb{R}^3}H\ast \eta^{k_n}\ast\rho^n_0(x)\rho^n_0(x)\,\dd x\right)\,\dd t\\
&=\left(\int_0^\infty\phi(t)\,\dd t\right)\left(\int_{\mathbb{R}^6}\sqrt{1+|v|^2} f_0(x,v)\,\dd v\,\dd x+\frac{\sigma_E}{2}\int_{\mathbb{R}^3}H\ast \rho_0(x)\rho_0(x)\,\dd x\right).
\end{split}\]
The case $\sigma_E=-1$ is subtler: by \eqref{liminfrelativistic} and \eqref{conservation1} we have that
\[\begin{split}
&\int_0^\infty\phi(t)\left(\int_{\mathbb{R}^6}\sqrt{1+|v|^2} f_t(x,v)\,\dd v\,\dd x-\frac{1}{2}\int_{\mathbb{R}^3}H\ast \rho_t(x)\rho_t(x)\,\dd x\right)\,\dd t\\
&\leq \liminf_{n\rightarrow \infty} \int_0^\infty \phi(t)\left(\int_{\mathbb{R}^6} \sqrt{1+|v|^2}f^n_t(x,v)\,\dd v\,\dd x-\frac{1}{2}\int_{\mathbb{R}^3}H\ast \rho_t(x)\rho_t(x)\,\dd x\right)\,\dd t\\
&\leq \left(\int_0^\infty\phi(t)\,\dd t\right)\left(\int_{\mathbb{R}^6}\sqrt{1+|v|^2} f_0(x,v)\,\dd v\,\dd x-\frac{1}{2}\int_{\mathbb{R}^3}H\ast \rho_0(x)\rho_0(x)\,\dd x\right)\\
&+\frac{1}{2}\limsup_{n\rightarrow\infty}\int_0^\infty \phi(t)\left(\int_{\mathbb{R}^3}H\ast \eta^{k_n} \ast\rho^n_t(x)\rho^n_t(x)-H\ast \rho_t(x)\rho_t(x)\,\dd x\right)\,\dd t
\end{split}\]
Notice that by \eqref{boundenergyT}, \eqref{ftregularity} and \eqref{calderonsobolev}, we have for every $T>0$,
\[
\sup_{t\in[0,T]}\|\rho_t\|_{L^{6/5}(\mathbb{R}^3)}+\sup_{n\in\mathbb{N}}\sup_{t\in[0,T]}\|\rho^n_t\|_{L^{6/5}(\mathbb{R}^3)}<\infty.
\]
Thus, by \Cref{proven}, we obtain that the last term equals $0$. Since $\phi$ was arbitrary and since $f_0$ has every energy bounded, we conclude that $f_t$ has every energy bounded for $\mathcal{L}^1$-almost every $t\in (0,\infty)$.

\textbf{Step 2: Bound on the total energy for every time.} Notice that the relativistic and electric potential energy is lower semicontinuous with respect to the strong $L^1_{\loc}$ and weak* $\mathcal{M}_+$ convergences, respectively. Hence, by the continuity of $t\longmapsto f_t\in L^1(\mathbb{R}^6)$ and $t\longmapsto \rho_t\in \mathcal{M}_+(\mathbb{R}^3)$ for the $L^1_{\loc}$ and weak* $\mathcal{M}_+$ convergences, respectively, combined with Step 1, we have that for $t_n\longrightarrow \bar{t}\in [0,\infty)$ such that \eqref{protobound} holds for all $t_n$, we may pass the limit and obtain \eqref{protobound} for $t=\bar{t}$.

\textbf{Step 3: Strong $\boldsymbol{L^1_{\loc}}$-continuity of the $\boldsymbol{\rho,\,J,\, E,\, B}$.} Given $t\in[0,\infty)$, let $t_n\longrightarrow t$. Fix $r>0$, and for any $R>0$
\[
\int_{B_r}\int_{\mathbb{R}^3}|f_{t_n}-f_t|\,\dd v\, \dd x\leq \int_{B_r}\int_{B_R}|f_{t_n}-f_t|\,\dd v\, \dd x+ R^{-1}\int_{B_r}\int_{\mathbb{R}^3}\sqrt{1+|v|^2}(f_{t_n}+f_t)\,\dd v\, \dd x.
\]
By the uniform boundedness of the relativistic energy with respect to time and the $L^1_{\loc}$ continuity of $f_t$, by taking the limit in $n$ and then in $R$, we conclude that $\rho_{t_n}\longrightarrow \rho_t$ in $L^1_{\loc}$. Moreover, since $|\hat{v}|<1$, we have
\[
\int_{B_r}|J_{t_n}-J_t|\, \dd x< \int_{B_r}\int_{\mathbb{R}^3}|f_{t_n}-f_t|\,\dd v\, \dd x\longrightarrow 0,
\]
thus $J_{t_n}\longrightarrow J_t$ in $L^1_{\loc}$. Finally, since $K\in L^1_{\loc}$ and $|J|(\mathbb{R}^3)<\rho(\mathbb{R}^3)<\infty$, we conclude that $E_t,\, B_t$ are also strongly continuous in $L^1_{\loc}(\mathbb{R}^3)$.

\textbf{Step 4: Globally defined flow.} We can combine the fact that $f_t$  has every energy bounded and \Cref{magneticpositive} to obtain that $E_t,\, B_t\in L^\infty([0,\infty);L^2(\mathbb{R}^3))$, thus by \Cref{maincorollary} we conclude that the trajectories of the maximal regular flow starting at any given $t$ do not blow up for $f_t$-almost every $(x,v)\in\mathbb{R}^6$.

\textbf{Step 5: Strong $\boldsymbol{L^1}$-continuity of $\boldsymbol{f}$.} By \Cref{existencesolution} and $L^1_{\loc}$-continuity of $f_t$, we deduce that finite energy solutions conserve mass, i.e., $\rho_t(\mathbb{R}^3)=\rho_0(\mathbb{R}^3)$ for every $t\in [0,\infty)$. In particular, solutions are strongly continuous in $L^1(\mathbb{R}^6)$ and not only $L^1_{\loc}(\mathbb{R}^6)$ (see \cite[Theorem 4.10]{vlasovpoisson}).
\endproof

\appendix
\section{Derivation}\label{deriv-model}

The relativistic Vlasov equation describes the evolution of a function $f:(0,\infty)\times \mathbb{R}^3\times \mathbb{R}^3\longrightarrow [0,\infty)$ under the action of a self-consistent acceleration $A:(0,\infty)\times \mathbb{R}^3\times \mathbb{R}^3\longrightarrow \mathbb{R}^3$:
\begin{equation}\label{vlasov}
	\partial_t f_t(x,v) +\hat{v}\cdot \nabla_x f_t(x,v)+A_t(x,v)\cdot \nabla_v f_t(x,v)=0 \quad \text{ in }\quad (0,\infty)\times \mathbb{R}^3\times \mathbb{R}^3.
\end{equation}
In this paper, we consider the acceleration given by
\[
A_t(x,v)=g_t(x)+ \frac{q}{m} (E_t(x)+\hat{v}\times B_t(x)),
\]
where $g_t$, $E_t$, and $B_t$ are the Newtonian gravitational, electric, and magnetic fields, respectively, and $q$ and $m$ are the particle charge and mass. Newtonian gravity implies that $g_t=Gm\nabla(-\Delta)^{-1}\rho_t$, where $G$ is the gravitational constant and $\rho_t$ the density of particles. We study the case in which the electromagnetic field satisfies one of the quasi-static limits of Maxwell's equations (see, for instance, \cite{manfredi} and references therein):
\begin{equation}\label{QES}
	\nabla\cdot E_t=\frac{q}{\epsilon_0}\rho_t,\quad \nabla \cdot B=0,\quad \nabla\times E_t=0,\quad \nabla\times B_t= \frac{q}{\epsilon_0}J_t+\partial_t E,
\end{equation}
or 
\begin{equation}\label{QMS}
	\nabla\cdot E_t=\frac{q}{\epsilon_0}\rho_t,\quad \nabla \cdot B=0,\quad \nabla\times E_t=-\partial_t B_t \quad \nabla\times B_t=\frac{q}{\epsilon_0} J_t,
\end{equation}
where $J_t$ is the relativistic particle current density. Equations \eqref{QES} and \eqref{QMS} are known as the quasi-electrostatic (QES) and quasi-magnetostatic (QMS) limits, respectively. The solution of \eqref{QES} can be written as
\[
E_t= -\frac{q}{\epsilon_0}\nabla(-\Delta)^{-1}\rho_t  \qquad \text{and} \qquad B_t=\frac{q}{\epsilon_0}\nabla\times(-\Delta)^{-1}J_t,
\]
while the solution of \eqref{QMS} is 
\[
E_t= -\frac{q}{\epsilon_0}\nabla(-\Delta)^{-1}\rho_t-\frac{q}{\epsilon_0}\partial_t (-\Delta)^{-1}J_t \qquad \text{and} \qquad B_t=\frac{q}{\epsilon_0}\nabla\times(-\Delta)^{-1}J_t.
\]
Notice that the leading term in the QES limit is the electric field whereas in the QMS it is the magnetic field. Hence, in the QES case, we can write $A_t$ in terms of $\rho_t$ and $J_t$ only:
\[
A_t(x,v)=\left(\frac{q^2}{4\pi \,\epsilon_0 m}-Gm\right)\int_{\mathbb{R}^3}\rho_t(y)\frac{x-y}{|x-y|^3}\,\dd y + \frac{q^2}{4\pi\, \epsilon_0 m} \, \hat{v} \times \int_{\mathbb{R}^3}J_t(y)\times \frac{x-y}{|x-y|^3}\,\dd y,
\]
where $\epsilon_0$ is the electric permittivity. Next, define the critical charge $q_c$ as
\[
q_c\coloneqq \pm \sqrt{4\pi \epsilon_0 G}\,m.
\] 
If $q>q_c$, we have that the electric field is stronger and, up to redefining of $\rho_t$ and $J_t$, we may write the acceleration as
\[
A_t(x,v)= \int_{\mathbb{R}^3}\rho_t(y)K(x-y)\,\dd y+\hat{v}\times\int_{\mathbb{R}^3}J_t(y)\times K(x-y)\,\dd y.
\] Analogously, if $q<q_c$, we can write
\[
A_t(x,v)= -\int_{\mathbb{R}^3}\rho_t(y)K(x-y)\,\dd y+\hat{v}\times\int_{\mathbb{R}^3}J_t(y)\times K(x-y)\,\dd y.
\]
In both cases, if we drop the magnetic field (since it is a lower order term), we have the relativistic Vlasov-Poisson system. Moreover, notice that in the critical case $q=q_c$, we only have the magnetic force acting in the evolution equation \eqref{vlasov}, which is exactly the same as if we only consider the leading term in the QMS limit, that is, the relativistic Vlasov-Biot-Savart system.

\end{document}